\documentclass[12pt,leqno]{article}
\usepackage{amsmath, amsthm, amsfonts, amssymb, color}
\usepackage{mathrsfs}
\newtheorem{thm}{Theorem}[section]

\newtheorem{lem}[thm]{Lemma}

\newtheorem{Remark}[thm]{Remark}
\theoremstyle{definition}

\newcommand{\scr}[1]{\mathscr #1}
\definecolor{wco}{rgb}{0.5,0.2,0.3}

\numberwithin{equation}{section}

\newcommand{\ua}{\uparrow}

\title{{\bf Stochastic Heat Equations with Values in a Manifold  via Dirichlet Forms}\footnote{Supported in
 part by  NSFC (11771037, 11671035, 11371099). Financial support by the DFG through the CRC 1283 ``Taming uncertainty and profiting from randomness and low regularity in analysis, stochastics and their applications'' is acknowledged.
}}
\author{
{\bf    Michael R\"{o}ckner$^{a)}$, Bo Wu$^{b,e)}$, Rongchan Zhu$^{a,c)}$\thanks{Corresponding author}, Xiangchan Zhu$^{a,d,}$}
\thanks{E-mail address: roeckner@math.uni-bielefeld.de(M.R\"{o}ckner), wubo@fudan.edu.cn(B.Wu),
zhurongchan@126.com(R.C.Zhu), zhuxiangchan@126.com(X.C.Zhu)}
\\
\footnotesize{  $^{a)}$ Department of Mathematics, University of Bielefeld, D-33615 Bielefeld, Germany}\\
\footnotesize {$^{b)}$ School  of Mathematical Sciences, Fudan
University, Shanghai 200433, China}\\
 \footnotesize{ $^{c)}$Department of Mathematics, Beijing Institute of Technology, Beijing 100081,  China}\\
\footnotesize{ $^{d)}$School of Science, Beijing Jiaotong University, Beijing 100044, China}
\\
 \footnotesize{ $^{e)}$ Institute for Applied Mathematics, University of Bonn,
Bonn 53115, Germany
}}

\date{}

\begin{document}
\maketitle

\def\paral{/\kern-0.55ex/}
\def\parals_#1{/\kern-0.55ex/_{\!#1}}
\def\R{\mathbb R} \def\EE{\mathbb E} \def\P{\mathbb P}\def\Z{\mathbb Z} \def\ff{\frac} \def\ss{\sqrt}
\def\H{\mathbb H}
\def\HH{\mathbf{H}}
\def\DD{\Delta} \def\vv{\varepsilon} \def\rr{\rho}
\def\<{\langle} \def\>{\rangle} \def\GG{\Gamma} \def\gg{\gamma}
\def\ll{\lambda} \def\LL{\Lambda} \def\nn{\nabla} \def\pp{\partial}
\def\dd{\text{\rm{d}}}
\def\Id{\text{\rm{Id}}}\def\loc{\text{\rm{loc}}} \def\bb{\beta} \def\aa{\alpha} \def\D{\scr D}
\def\E{\scr E} \def\si{\sigma} \def\ess{\text{\rm{ess}}}
\def\beg{\begin} \def\beq{\beg}  \def\F{\scr F}
\def\Ric{\text{\rm{Ric}}}
\def\Vol{\text{\rm{Vol}}}
\def\Var{\text{\rm{Var}}}
\def\Ent{\text{\rm{Ent}}}
\def\Scal{\text{\rm{Scal}}}
\def\Hess{\text{\rm{Hess}}}\def\B{\scr B}
\def\e{\text{\rm{e}}} \def\ua{\underline a} \def\OO{\Omega} \def\b{\mathbf b}
\def\oo{\omega}     \def\tt{\tilde} \def\Ric{\text{\rm{Ric}}}
\def\cut{\text{\rm{cut}}} \def\P{\mathbb P} \def\K{\mathbb K}
\def\ifn{I_n(f^{\bigotimes n})}
\def\fff{f(x_1)\dots f(x_n)} \def\ifm{I_m(g^{\bigotimes m})} \def\ee{\varepsilon}
\def\C{\scr C}
\def\PP{\scr P}
\def\M{\scr M}\def\ll{\lambda}
\def\X{\scr X}
\def\T{\scr T}
\def\A{\mathbf A}
\def\LL{\scr L}\def\LLL{\Lambda}
\def\gap{\mathbf{gap}}
\def\div{\text{\rm div}}
\def\Lip{\text{\rm Lip}}
\def\dist{\text{\rm dist}}
\def\cut{\text{\rm cut}}
\def\supp{\text{\rm supp}}
\def\Cov{\text{\rm Cov}}
\def\Dom{\text{\rm Dom}}
\def\Cap{\text{\rm Cap}}\def\II{{\mathbb I}}\def\beq{\beg{equation}}
\def\sect{\text{\rm sect}}\def\H{\mathbb H}

\begin{abstract}In this paper, we prove the existence of martingale solutions to the stochastic heat equation taking values in a
Riemannian manifold, which admits Wiener
(Brownian bridge) measure on the Riemannian path (loop) space as an invariant measure using a suitable Dirichlet form. Using the Andersson-Driver approximation, we heuristically derive a form of the equation solved by the process given by the Dirichlet form.

Moreover, we  establish
the log-Sobolev inequality for the Dirichlet form in the path space. In addition, some characterizations for the lower  bound of the Ricci curvature are presented related to the stochastic heat equation.
\end{abstract}

\noindent Keywords: Stochastic heat equation; Ricci curvature; Functional inequality; Quasi-regular Dirichlet form
\vskip 2cm

\section{Introduction}\label{sect1}
This work is motivated by Tadahisa Funaki's pioneering work \cite{Fun92} and Martin Hairer's recent work \cite{Hai16}. The former had proved the existence and uniqueness of a
natural evolution driven by regular noise  on loop space over a Riemannian manifold $(M,g)$ by the classical theory of stochastic differential equation, and the latter considered the singular noise case, i.e. the associated stochastic heat equation may be interpreted formally as
\begin{equation}\dot{u}=-\nabla E (u)+\sum^m_{i=1}\sigma_i(u)\xi_i\end{equation}
for solutions $u(t,\cdot):S^1\rightarrow M$, i.e.
the formal Langevin dynamics
associated to the energy given by
$$E(u)=\int_{S^1} g_{u(s)}(\partial_su(s),\partial_su(s))\dd s,$$
and $\sum^m_{i=1}\sigma_i(u)\xi_i$ is a suitable ``white noise on loop space".
By Andersson-Driver's work in \cite{AD99}, we know that there exists an explicit relation between the Langevin energy $E(u)$ and the Wiener (Brownian bridge) measure (see also  \cite{Li,TL, W2}). In \cite{AD99}, the Wiener
(Brownian bridge) measure $\mu$ has been interpreted as the limit of a natural approximation of the  measure
$\exp(-\frac{1}{2}E(u))\D u$, where $\D u$ denotes a `Lebesgue' like measure on path space.
By observing the above connection, one may think the solution to the stochastic heat equation (1.1) may have $\mu$ as an invariant (even symmetrizing) measure.
In this paper, starting from the invariant measure $\mu$ we use the theory of Dirichlet forms to construct a natural evolution which admits $\mu$ as an invariant measure.

Actually, on a Riemannian path/loop space there exists  another process which also admits $\mu$ as an invariant measure associated with the Dirichlet form $\E^{OU}$ given by the Malliavin gradient, which sometimes is called the \textbf{O-U Dirichlet process}. These processes had first been constructed by Driver and R\"{o}ckner in \cite{DR92} for pinned loop space and by  Albeverio,  L\'{e}andre and  R\"{o}ckner in \cite{ALR93} for the free loop space by using Dirichlet form theory. After that there were several follow-up papers concentrating on more general cases,
see \cite{EM, L, WW08, CW14}. In addition, Norris in \cite{Nor95,Nor98} obtained some similar processes based on the theory of stochastic differential equations.

Recently, Hairer \cite{Hai16} wrote equation (1.1) in  local coordinates heuristically as
\begin{equation}\label{eq1.1}\partial_t{u}^\alpha=\partial_s^2u^\alpha+\Gamma_{\beta\gamma}^\alpha(u)\partial_su^\beta\partial_su^\gamma+\sigma_i^\alpha(u)\xi_i,\end{equation}
where Einstein summation convention over repeated indices is applied and
$\Gamma_{\beta\gamma}^\alpha$ are the Christoffel symbols for the Levi-Civita connection of $(M, g)$. $\sigma_i$ are some suitable vector fields on $M$. This equation may be considered as a certain kind of multi-component version of the KPZ equation.
By the theory of regularity structures recently developed in \cite{Hai14, BHZ16, CH16}, local well-posedness of \eqref{eq1.1} has been obtained in \cite{Hai16}.

 When $M=\mathbb{R}^d$, the process constructed by Driver and R\"{o}ckner in \cite{DR92} is the O-U process in the Mallivan calculus, whereas equation \eqref{eq1.1} corresponds to the stochastic heat equation.
To explain the difference of the above two processes, let us first consider the following two stochastic equations on $L^2([0,1];\mathbb{R}^d)$:

\noindent\textbf{A. Ornstein-Unlenbeck process}
$$\dd X(t)=\frac{1}{2}X(t)\dd t+(-\Delta_D)^{-\frac{1}{2}}\dd W(t),$$
\textbf{B. Stochastic heat equation}
$$\dd X(t)=\frac{1}{2}\Delta_{D}X(t)\dd t+\dd W(t),$$
where $\Delta_{D}$ is the operator $\frac{\dd^2}{\dd s^2}$ on $L^2([0,1];\mathbb{R}^d)$ with boundary condition $h(0)=0, h(1)=0$ and $W$ is an $L^2$-cylindrical Wiener process.
Solutions to these two equations share the same Gaussian invariant measure $N(0,(-\Delta_{D})^{-1})$ in $L^2([0,1];\mathbb{R}^d)$. It is well known that $N(0,(-\Delta_{D})^{-1})$ has full (topological) support on $C([0,1];\mathbb{R}^d)$ and is the same as the law of the Brownian bridge on $C([0,1];\mathbb{R}^d)$ starting from $0$.

In the first part of this paper, we construct the solutions to the stochastic heat equation taking values in  a Riemannian manifold $M$ by Dirichlet form theory (i.e. we replace in Case \textbf{B} above $\mathbb{R}^d$ by $M$). Unlike as in \cite{DR92}, we consider the closure $(\E,\D(\E))$ of the following bilinear form with the  reference measure $\mu=$ the law of Brownian motion on $M$ (path space case)/the law of Brownian bridge on $M$ (loop space case):
$$\E(F,G):=\frac{1}{2}\int_{\mathbf{E}} \langle DF, DG\rangle_{\mathbf{H}}\dd\mu=\frac{1}{2}\sum_{k=1}^\infty\int_{\mathbf{E}} D_{h_k}F D_{h_k}G  \dd\mu; \quad F,G\in\F C_b^1,$$
where $\F C_b^1$ is introduced in (2.2) below, ${\mathbf{H}}:=L^2([0,1];\mathbb{R}^d)$, $\mathbf{E}$ is introduced in Section 2.1 and $DF$ is the $L^2$-derivative defined in Section 2 with $\{h_k\}$ being an orthonormal basis in ${\mathbf{H}}$. In
this case, we call the associated Dirichlet form \textbf{{$\mathbf{L^2}$}-Dirichlet form}. When $M=\mathbb{R}^d$,
this Dirichlet form just corresponds to the stochastic heat equation, i.e. Case \textbf{B} above, while in \cite{DR92} the considered Dirichlet form corresponds to the process from Case A above, if $M=\mathbb{R}^d$. Therefore, below we denote the Dirichlet form in \cite{DR92} by $\E^{OU}$.  By simple computations, one sees that the classical cylinder test functions $u(\gamma)=f(\gamma_t)$ considered in \cite{DR92} are not in the domain of $\D(\E)$, since $\E(u,u)$ might be infinity. Thus, we need to choose a class of suitable functions $\F C_b^1$ introduced in \eqref{eq2.2} below. As usual in Dirichlet from theory to obtain the corresponding Markov process, i.e. in our case the solution to the stochastic heat equation in the path/loop space over $M$, we have to prove \textbf{a) closability }and \textbf{b) quasi-regularity}, which is done in Section 2 below (see Theorem 2.2).  To prove the quasi-regularity of the closure of $\E$, the uniform distance will be replaced by $L^2$-distance mentioned in Subsection 2.1. Then we  obtain martingale solutions to the stochastic heat equations, which admit $\mu$ as an invariant measure on path space and loop space, respectively.

We would like to emphasize that both the Dirichlet form approach and regularity structure theory can be used to construct rigorously the natural Markov process associated with \eqref{eq1.1} (see Section 4).  The Markov process constructed by the Dirichlet form approach is a global martingale solution starting from quasi-every starting point (which is a path/loop), which has the law of the Brownian motion resp. Brownian bridge on $M$ as an invariant measure. The solution obtained by regularity structure in \cite{Hai16} is a local strong solution to \eqref{eq1.1} starting from every point ($=$path/loop).

In this paper we consider four different cases: pinned path resp. loop spaces and free path resp. loop spaces.
For a better understanding of the measure and the stochastic heat equation on these spaces, let us first look at the simplest case : $M=\mathbb{R}^d$.
For the case of the path space, the reference measure is  $P_x:= N(x,(-\Delta_{D,N})^{-1})$, which is the law of Brownian motion starting from a fixed point $x\in\mathbb{R}^d$,  where $\Delta_{D,N}$ is the operator $\frac{\dd^2}{\dd s^2}$ on $L^2([0,1];\mathbb{R}^d)$ with boundary condition $h(0)=0, h'(1)=0$,
and the corresponding SPDE constructed by the $L^2$-Dirichlet form $\E$ above is the following:
\begin{equation}\label{eqflat}\dd X(t)=\frac{1}{2}\Delta_{D,N}(X(t)-x)\dd t+\dd W(t),\end{equation}
with $W$ an $L^2([0,1];\mathbb{R}^d)$-cylindrical Wiener process.
On the loop space, the reference measure is $P_{x,x}= N(x,(-\Delta_{D})^{-1})$, which is the law of Brownian bridge starting from $x\in\mathbb{R}^d$, then
the corresponding SPDE is:
$$\dd X(t)=\frac{1}{2}\Delta_{D}(X(t)-x)\dd t+\dd W(t),$$ where $\Delta_{D}$ is the operator $\frac{\dd^2}{\dd s^2}$ on $L^2([0,1];\mathbb{R}^d)$ with boundary condition $h(0)=0, h(1)=0$.
 For the case of free path/loop space, we have  the following:  Let $\sigma$ be a probability measure on $\R^d$. Then  the reference measure for the free path case is given  by $\int P_x\sigma(\dd x)= \int N(x,(-\Delta_{D,N})^{-1})\sigma(\dd x)$ and
the corresponding SPDE is:
$$\dd X(t,x)=\frac{1}{2}\Delta_{D,N}(X(t,x)-X(t,0))\dd t+\dd W(t).$$
Similarly, the reference measure for the free loop case is given   by $\int P_{x,x}\sigma(\dd x)= \int N(x,(-\Delta_{D})^{-1})\sigma(\dd x)$ and
the corresponding SPDE is:
$$\dd X(t,x)=\frac{1}{2}\Delta_{D}(X(t,x)-X(t,0))\dd t+\dd W(t).$$

In the second part of this paper, we use functional inequalities  to study the properties of the solutions to the stochastic heat equations on path space over a Riemannian manifold $M$.   Functional inequalities for Ornstein-Unlenbeck process on Riemannian path space have been well-studied (see
\cite{F94, AE95, ALR93, EL99, AE95, H97, H99, W04, CW14} and references therein).
Since the $L^2$-Dirichlet form associated with the stochastic heat equation is larger than the O-U Dirichlet form $\E ^{OU}$ constructed in \cite{DR92} (i.e., $\E ^{OU}(u,u)\leq \E (u,u)$ for $u\in D(\E)$),
all the functional inequalities with respect to $\E ^{OU}$ still hold in the stochastic heat equation case, which implies that the former requires stronger Ricci curvature conditions than the latter. In fact,
from recent results in \cite{N} by Naber, we know that the Poincare inequality/log-Sobolev inequality for the twisted O-U Dirichlet form requires the uniformly bounded Ricci curvature. And for the $L^2$-Dirichlet form, it only needs lower bounded  Ricci curvature, which had already been proved before by Gourcy-Wu in \cite{GW06}. In this section we also establish the log-Sobolev inequality for $\E$, but our constant $C(K)$ is smaller than Gourcy-Wu's constant (see Theorem \ref{T3.1} below). In particular, when $M$ is an Einstein manifold with constant Ricci curvature $K\in\R$, our constant $C(K)$ in the log-Sobolev inequality is optimal in the sense that $\lim_{K\rightarrow0}C(K)=\frac{4}{\pi^2}$ and $\frac{4}{\pi^2}$ is the optimal constant for the log-Sobolev inequality in the flat case (see Theorem 3.3 below). Here we want to emphasize that the log-Sobolev inequality implies the $L^2$-ergodicity of the solution to the stochastic heat equation (see Remark 3.2 below for other consequences).

As mentioned above, the log-Sobolev inequality is a consequence of a geometric property of the manifold. It is very interesting to ask to what extent these geometric properties are also necessary for the log-Sobolev inequality to hold for $\E$ as above. The most interesting work is related to the Bakry-Emery criterion, which gives a characterization of the lower boundedness of the Ricci curvature in terms of the log-Sobolev inequality for the classical Dirichlet form on a Riemannian manifold (see \cite{BE85}). Recently, Naber in \cite{N} characterizes  uniform boundedness  of the Ricci curvature  using the O-U process on path space. Wang-Wu\cite{WW16} obtained a more general characterization of the Ricci curvature and the second fundamental form on the boundary of the Riemannian manifold  using a new method. After that, this result has been extended to general uniform bounds of the Ricci curvature by Wu \cite{W1} and Cheng-Thalimaier \cite{CT}. In addition, Wu \cite{W1} and Wang \cite{W17} gave some characterization for the upper bound of the Ricci curvature by analysis on path space and the Weitzenb\"{o}ck-Bochner integration formula respectively.
Similar to the above case,
in Subsection 3.2 we  give
some (equivalent) characterizations of the lower boundedness  of the Ricci curvature by using the $L^2$-gradient and a properly weighted decomposition of the $L^2$-Dirichlet form on path space.

In the last part of the paper, we discuss the form of the stochastic heat equations constructed in Section 2.
Dirichlet form theory is a useful tool to construct  stochastic processes on infinite dimensional spaces (see \cite{AR91, MR92}). In the flat case  we can  use Dirichlet form theory to write an SPDE which the process satisfies (see \cite{AR91} for the $p(\Phi)_2$ model). This helps  to  obtain
new properties of the $\Phi^4_2$
field (see \cite{RZZ15, RZZ16}).   However, in the Riemannian manifold case, the explicit form of the SPDE cannot be deduced directly since there are no linear functions on the Riemannian manifold. It will  be seen from Section 2 that the martingale part is  space-time white noise and thus is very rough. To define the drift part,  renormalization is required (see \cite{Hai14}). In Section 4 we construct suitable approximation processes on the piecewise geodesic space  using the approximation measures from \cite{AD99} and discuss the convergence of the approximations, which gives two limiting   forms of the stochastic heat equations. One is more related to the Mosco convergence of the approximating Dirichlet forms to the $L^2$-Dirichlet form $\E$ and the integration by parts formula obtained by Driver (see \eqref{eq4.37} and Remark \ref{r4.9} below). The other is more related to \eqref{eq1.1} constructed by Hairer ( see \eqref{eq} below).  The different forms of the equations originate from the different choices of the diffusion coefficients. We conjecture that the Markov process constructed in Section 2 by the Dirichlet form $\E$ has the same law as the solution to \eqref{eq1.1} constructed by Hairer in \cite{Hai16}. To obtain the Mosco convergence of the corresponding Dirichlet forms  in general Markov uniqueness of the limiting $L^2$-Dirichlet form is required, which is a difficult problem in Dirichlet form theory.   We hope to be able to use the theory of regularity structures/ paracontrolled distribution method to  make the heuristic convergence of the corresponding solutions in Section 4 rigorous in our future work.

This paper is organized as follows: In Section 2 we construct the $L^2$-Dirichlet form $\E$ on the pinned (free) path/loop space. By this we obtain existence of martingale solutions to the stochastic heat equation on path space and loop space. In Section 3.1, we derive functional inequalities for the $L^2$-Dirichlet form $\E$.
The equivalent characterizations of the lower boundedness  of the Ricci curvature are obtained in Section 3.2. In Section 4.2, we construct approximation processes on the piecewise geodesic space by considering Dirichlet forms with respect to the approximation measure.  In Section 4.3 we discuss the convergence of the approximation processes and  the form of the limiting stochastic equation heuristically.

\section{Construction of Dirichlet form}\label{sect2}
\subsection{Dirichlet form on pinned path Space}
Throughout this article, suppose that $M$ is a complete and stochastic complete Riemannian manifold with dimension $d$, and $\rho$ is the Riemannian distance on $M$. In this section we assume that $M$ is compact for simplicity and for the more general case, we refer to Remark 2.1.
Fix $o\in M$ and $T>0$,
the based path space $W^T_{o}(M)$ over $M$ is defined by
$$W^T_{o}(M):=\{\gamma\in C([0,T];M):\gamma(0)=o\}.$$
Then $W^T_{o}(M)$ is a Polish
space under the uniform distance
$$d_\infty(\gamma,\sigma):=\displaystyle\sup_{t\in[0,T]}\rho(\gamma(t),\sigma(t)),\quad\gamma,\sigma\in W^T_{o}(M).$$
For convenience, we write $W_{o}(M):=W^1_{o}(M)$. In the following we consider $W_o(M)$ for simplicity.
In order to construct Dirichlet forms associated to stochastic heat equations on Riemannian path space, we first need to introduce the following $L^2$-distance, which is a smaller distance than the above uniform distance $d_\infty$ on  $W_{o}(M)$:
$$\tilde{d}(\gamma,\eta)^2:={\int_0^1\rho(\gamma_s,\eta_s)^2\dd s}, \quad \gamma,\eta\in W_o(M).$$
The $L^2$-distance $\tilde{d}$ is quite crucial to prove the quasi-regularity for the  Dirichlet form mentioned in Theorem \ref{T2.1}.
Let ${\bf E}$ be the closure of $W_o(M)$ in $$\left\{\eta:[0,T]\rightarrow M; \int_0^1\rho(o,\eta_s)^2\dd s<\infty\right\}$$ with respect to the distance $\tilde{d}$, then ${\bf E}$ is a Polish space.

Before stating our main results in this section, let us  recall some basic notation and introduce the Brownian motion on $M$. Let $\nabla$ be the Riemannian connection on $M$ and the curvature tensor $R$ of $\nabla$ is given by
$$R(X,Y)Z=\nabla_X\nabla_YZ-\nabla_X\nabla_YZ-\nabla_{[X,Y]}Z$$
for all vector fields $X, Y$ and $Z$ on $M$. The Ricci tensor $\Ric$ and  the scalar curvature $\Scal$ of $M$ are traces of $R$  and $\Ric$ respectively, i.e.,
$$\Ric X:=\sum_{i=1}^d R(X,\bar{e}_i)\bar{e}_i, \quad\textrm{Scal}=\sum_{i=1}^d \langle\Ric\bar{e}_i,\bar{e}_i\rangle,$$
where $\{\bar{e}_i\}$ is an orthonormal frame.

Let $O(M)$ be the orthonormal frame bundle over $M$, and let $\pi: O(M) \rightarrow M$ be the canonical projection.
Furthermore, we choose a standard othornormal basis $\{H_i\}_{i=1}^d$
of horizontal vector fields on $O(M)$ and consider the following SDE:
\begin{equation}\label{eq2.1}
\begin{cases}
&\dd U_t=\sum^d_{i=1}H_i(U_t)\circ \dd B_t^i,\quad t\geq0\\
& U_0=u_o,
\end{cases}
\end{equation}
where $u_o$ is a fixed orthonormal basis of $T_o M$ and $B^1_t,\cdots,B_t^d$ are independent
Brownian motions on $\mathbb{R}$.
Then  $x_t:=\pi(U_t),\ t\geq0,$ is the Brownian motion on
$M$ with initial point $o$, and $U_{\cdot}$ is the (stochastic) horizontal lift along
$x_{\cdot}$. Let $\mu$ be the distribution of $x_{[0,1]}:=\{x(t)|t\in[0,1]\}$. Then $\mu$
is a probability measure on $W_{o}(M)$.

In the following we use $\langle\cdot,\cdot\rangle$ to denote the inner product in $\mathbb{R}^d$.

Let $\F C_b^1$ be a space of
$C_b^1$ cylinder functions on ${\bf E},$ defined as follows: for every $F\in \F C_b^1$, there exist some $m\geq1, ~m\in \mathbb{N},~ f\in C_b^1(\mathbb{R}^m), g_i\in C_b^{0,1}([0,1]\times M)$, $i=1,...,m$,
such that
\begin{equation}\label{eq2.2}\aligned
F(\gamma)=f\left(\int_0^1  g_1(s,\gamma_s) \dd s,\int_0^1  g_2(s,\gamma_s) \dd s,...,\int_0^1  g_m(s,\gamma_s) \dd s\right),\quad \gamma\in {\bf E}.\endaligned\end{equation}
Here $C_b^{0,1}([0,1]\times M)$ denotes the functions which are continuous w.r.t. the first variable and $C^1$-
differentiable w.r.t. the second variable with continuous derivatives. It is easy to see that $\F C_b^1$ is dense in $L^2({\bf E}, \mu)$.
For any $F\in \F C_b^1$ of the form \eqref{eq2.2} and $h\in \HH:=L^2([0,1];\mathbb{R}^d)$, the directional derivative of $F$ with respect to $h$ is given by
$$D_hF(\gamma)=\sum_{j=1}^m\hat{\partial}_jf(\gamma)\int_0^1
 \left\langle U_s^{-1}(\gamma)\nabla g_j(s,\gamma_s),h_s
 \right\rangle \dd s,\quad \gamma\in W_o(M),$$
where
$$\hat{\partial}_jf(\gamma):=\partial_jf\bigg(\int_0^1  g_1(s,\gamma_s) \dd s,\int_0^1  g_2(s,\gamma_s) \dd s,...,\int_0^1  g_m(s,\gamma_s) \dd s\bigg).$$
and $\nabla g_j$ denotes the gradient w.r.t. the second variable.
Without loss of generality, for $\gamma\in \mathbf{E}\backslash W_o(M)$ we  take $D_hF(\gamma)=0$. By the Riesz representation theorem,
there exists a gradient operator $DF(\gamma)\in \HH$ such that $\langle DF(\gamma),h\rangle_{\HH}=D_hF(\gamma), \gamma\in \mathbf{E}, h\in \HH$. In particular, for $\gamma\in W_o(M)$,
\begin{equation}\label{eq2.3}\aligned DF(\gamma)(s)=\sum_{j=1}^m\hat{\partial}_jf(\gamma)U_s^{-1}(\gamma)\nabla g_j(s,\gamma_s).\endaligned\end{equation}

\begin{Remark} 
In fact, for a more general Riemannian manifold the main results in this section still hold. But when we prove the quasi-regularity of $\E$, it is  required that the function $g$ is allowed to be the distance function, which can be approximated by $C^1$-functions in a suitable way.  This will be considered in a forthcoming paper.
\end{Remark}

Let $\H$ denote the Cameron-Martin space:
$$\mathbb{H}:=\left\{h\in C^1([0,1];\mathbb{R}^d)\Big| h(0)=0, \|h\|^2_\H:=\int^1_0|h'(s)|^2\dd s<\infty\right\}.$$
Taking a sequence of elements $\{h_k\}\subset \H$ such that it is an orthonormal basis in $\HH$, consider the following symmetric quadratic form
$$\E(F,G):=\frac{1}{2}\int_{\bf E} \langle DF, DG\rangle_{\HH}\dd\mu=\frac{1}{2}\sum_{k=1}^\infty\int_{\bf E} D_{h_k}F D_{h_k}G  \dd\mu; \quad F,G\in\F C_b^1.$$
The following is the main results in this section.

\beg{thm}\label{T2.1}  The quadratic form $(\E, \F C_b^1)$
is closable and its closure $(\E,\D(\E))$ is a quasi-regular Dirichlet form on $L^2({\bf E};\mu)=L^2(W_o(M);\mu)$.
\end{thm}

The proof of Theorem \ref{T2.1} will be given at the end of this subsection.
Using the theory of Dirichlet forms (see \cite{MR92}), we obtain the following associated diffusion process.

 \beg{thm}\label{T2.2} There exists a conservative (Markov) diffusion process
 $M=(\Omega,\F,(\M_t),$ $(X_t)_{t\geq0},(\mathbf{P}^z)_{z\in \mathbf{E}})$ on ${\bf E}$ \emph{properly associated with} $(\E,\D(\E))$, i.e. for $u\in L^2({\bf E};\mu)\cap\B_b({\bf E})$, the transition semigroup $P_tu(z):={\bf E}^z[u(X_t)]$ is an $\E$-quasi-continuous version of $T_tu$ for all $t >0$, where $T_t$ is the semigroup associated with $(\E,\D(\E))$.
\end{thm}

Here $\B_b({\bf E})$ denotes the set of the bounded Borel-measurable functions and for the notion of $\E$-quasi-continuity we refer to \cite[ChapterIII, Definition 3.2]{MR92}. Moreover, by the Fukushima decomposition we have:
\beg{thm}\label{T2.4}  There exists a \emph{properly  $\E$-exceptional set} $S\subset {\bf E}$, i.e. $\mu(S)=0$ and $\mathbf{P}^z[X_t\in {\bf E}\setminus S, \forall t\geq0]=1$ for $z\in E\backslash S$, such that $\forall z\in {\bf E}\backslash S$ under $\mathbf{P}^z$,  the sample paths of the associated  process $M=(\Omega,\F,(\M_t),$ $(X_t)_{t\geq0},(\mathbf{P}^z)_{z\in {\bf E}})$ on ${\bf E}$ satisfy the following for $u\in \D(\E)$
 \begin{equation}\label{eq2.4}\aligned u(X_t)-u(X_0)=M_t^u+N_t^u\quad \mathbf{P}^z-a.s.,\endaligned\end{equation}
 where $M^u$ is a martingale with quadratic variation process given by $\int_0^t |Du(X_s)|_\HH^2\dd s$ and $N_t$ is a zero quadratic variation process. In particular, for $u\in D(L)$, $N_t^u=\int_0^tLu(X_s)ds$, where $L$ is the generator of $(\E,\D(\E))$.
\end{thm}

\beg{Remark}\label{r2.5}
i) If we choose $u(\gamma)=\int_{r_1}^{r_2}u^\alpha(\gamma_s)ds\in\F C_b^1$, $0\leq r_1<r_2\leq 1$, with $u^\alpha$ being  local coordinates on $M$, then the quadratic variation process for $M^u$ is the same as  that for the martingale part in (1.2) (see Remark \ref{r4.9}).

ii) Although the stochastic horizontal lift $(U_t(\gamma))_{t\in [0,\infty)}$ is applied
in the definition of $(\E, \F C_b^1)$, the value of $\E(F,F)$ is independent of $(U_t(\gamma))_{t\in [0,\infty)}$.
In particular, by the definition \eqref{eq2.3} of the gradient, we have
$$\E(F,G)=\frac{1}{2}\int
\sum_{i=1}^m\sum_{j=1}^l\hat{\partial}_if_1(\gamma)\hat{\partial}_jf_2(\gamma)\int^{1}_0\langle \nabla g^1_i(s,\gamma_s),\nabla g^2_j(s,\gamma_s)\rangle \dd s\dd\mu$$
for any $F,G\in\F C_b^1$ with
$$\aligned
&F(\gamma)=f_1\left(\int_0^{1}  g^1_1(s,\gamma(s)) \dd s,\int_0^{1}   g^1_2(s,\gamma(s)) \dd s,...,\int_0^{1}   g^1_m(s,\gamma(s)) \dd s\right)\\
&G(\gamma)=f_2\left(\int_0^{1}  g^2_1(s,\gamma(s)) \dd s,\int_0^{1}   g^2_2(s,\gamma(s)) \dd s,...,\int_0^{1}   g^2_l(s,\gamma(s)) \dd s\right),\quad \gamma\in {\bf E}.\endaligned$$
This implies the quadratic form $\E$ is independent of $(U_t(\gamma))_{t\in [0,\infty)}$. This is  different from the O-U Dirichlet form, since the latter depends on the parallel translation $(U_t(\gamma))_{t\in [0,\infty)}$.

iii) In Section 4 we have another reference measure $\mu_0:=e^{-\frac{1}{6}\int_0^1\textrm{\Scal}(\gamma(s))\dd s}\dd \mu(\gamma)$, which is also related to $\exp(-\frac{1}{2}E(u))\D u$ mentioned in introduction.  We can also construct the $L^2$-Dirichlet form $(\E^0,\D(\E^0))$ with respect to $\mu_0$ and the  results in Theorems 2.2-2.4 still hold in this case.
\end{Remark}

\ \newline\emph{\bf Proof of Theorem \ref{T2.1}.}
{\bf$(a)$ Closablity:}
By the integration by parts formula (refer to \cite{Dri92}, also see \cite{H95, H02}): for $h\in \mathbb{H}$,
\begin{equation}\label{eq2.5}\int D_hF\dd\mu=\int F\beta_h\dd\mu \end{equation}
for every cylinder function depending on finite times $F(\gamma)=f(\gamma_{t_1},...,\gamma_{t_m})$,  where $f\in C_b^1(M^m)$ and $t_i\in [0,1], i=1,...,m,$
 $$L^2({\bf E},\mu)\ni\beta_h:=\int^1_0\left\langle h_s'+\frac{1}{2}\Ric_{U_s}(h_s),\dd B_s\right\rangle,$$
 where
$$\langle\Ric_{U_s}(a_1),a_2\rangle:=\langle\Ric(U_sa_1), U_sa_2\rangle_{T_{\gamma_s}M},\quad a_1,a_2\in\mathbb{R}^d.$$
For each $F(\gamma)=f(\int_0^1g_1(s,\gamma_s)\dd s, \cdots,\int_0^1g_m(s,\gamma_s)\dd s)\in\F C_b^1$,  choose $$F_n=f\left(\frac{1}{n}\sum_{i=1}^n g_1(i/n,\gamma_{i/n}),\cdots, \frac{1}{n}\sum_{i=1}^n g_m(i/n,\gamma_{i/n})\right).$$
Then $F_n$ and $D_hF_n$ $L^2$-converge to $F$ and $D_hF$ respectively. Thus, we deduce that \eqref{eq2.5} holds for $F\in \F C_b^1$.

Since $\beta_h\in L^2({\bf E},\mu)$, it is standard to prove that $(\E,\F C^1_b)$ is closable (see \cite{DR92} or \cite{L,WW08,CW14}). For the completeness of the proof we write it in  detail.
Let $\{F_n\}_{n=1}^{\infty}\subseteq \F C_b^1$ be a sequence of cylinder functions with
\begin{equation}\label{eq2.6}
\begin{split}
\lim_{n \rightarrow \infty}\mu\left[ F_n^2\right]=0,\ \
\lim_{n,m \rightarrow \infty}\E\left(F_n-F_m,F_n-F_m\right)=0.
\end{split}
\end{equation}
Thus $\{D F_n\}_{n=1}^{\infty}$ is a Cauchy sequence in
$L^2\left({\bf E}\rightarrow \HH;\mu\right)$ for which there exists a limit $\Phi$. It suffices to prove that $\Phi=0$. By  \eqref{eq2.5},
for $G \in \F C^1_b$ and $k\geq1$, we have
\begin{equation}\label{eq2.7}
\begin{split}
&\mu\left[\langle D F_n, h_k\rangle_{\HH}G \right]\\
&=\mu\left[\langle D \left(F_nG\right), h_k\rangle_{\HH}\right]
-\mu\left[\langle D G, h_k\rangle_{\HH}F_n\right]\\
&=\mu\left[F_nG \int_0^1\left\langle h_k'(s)+\frac{1}{2}\Ric_{U_s}
h_k(s)
, \dd B_s\right\rangle\right]-\mu\left[\langle D G, h_k\rangle_{\HH}F_n\right].
\end{split}
\end{equation}
Since $G$ and $D G$ are bounded and $\int_0^1\left\langle h_k'(s)+\frac{1}{2}\Ric_{U_s}
h_k(s), \dd B_s\right\rangle\in L^2({\bf E};\mu)$, $F_n$ converges to $0$ in $L^2(\mu)$, we may take the limit
$n \rightarrow \infty$ under the integral in \eqref{eq2.7} and conclude
\begin{equation*}
\begin{split}
&\mu\left[\langle \Phi, h_k\rangle_{\HH}G\right]=0,\quad \forall\
G \in \F C^1_b,\ k\ge 1,
\end{split}
\end{equation*}
which implies that $\Phi=0$, a.s., and that
$(\E,\F C^1_b)$ is closable. By  standard methods, we show easily that its closure $(\E,\D(\E))$ is a Dirichlet form.

{\bf$(b)$  Quasi-Regularity:} By the Nash embedding theorem we may assume that $M$ is embedded isometrically into $\mathbb{R}^N$ for a large enough $N\in \R$:
$$\psi: p\mapsto\psi(p)=(\psi^1(p),...,\psi^N(p))\in \mathbb{R}^N.$$ Then the distance $\rho(p,q)$ is equivalent to $\rho_0(p,q):=|\psi(p)-\psi(q)|$ for $p,q\in M$ and $\psi$ is smooth on $M$, which implies that the two distances $\tilde{d}(\gamma,\eta)^2$ and $\bar{d}(\gamma,\eta)^2:=\sum_{i=1}^N\int_0^1(\psi^i(\gamma(s))-\psi^i(\eta(s)))^2\dd s$ on the path space $ {\bf E}$ are equivalent to each other.
Since ${\bf E}$ is separable we can choose a fixed countable dense set $\{\xi_m|m\in\mathbb{N}\}\subset W_o(M)$ in ${\bf E}$.
We first prove the tightness of the capacity for $(\E,\D(\E))$:
Let $\varphi\in C_b^\infty(\mathbb{R})$ be an increasing function  satisfying
$$\varphi(t)=t,\quad \forall~t\in[-1,1]~~\text{and}~~\|\varphi'\|_{\infty}\leq 1.$$
And for $m\in\mathbb{N}$, the function $v_m:{ \bf E}\rightarrow\mathbb{R}$ is given by
$$v_m(\gamma)=\varphi(\bar{d}(\gamma,\xi_m)^2),\quad\gamma\in {\bf E}.$$
Suppose we can show that
\begin{equation}\label{eq2.8}w_n:=\inf_{m\leq n}v_m, n\in \mathbb{N}, \textrm{ converges } \E -\text{quasi-uniformly to zero on} ~{\bf E},\end{equation}
then for every $k\in\mathbb{N}$ there exists a closed set $F_k$ such that $\textrm{Cap}(F_k^c)<\frac{1}{k}$ and $w_n\rightarrow0$ uniformly on $F_k$. For every $0<\epsilon<1$ there exists $n\in\mathbb{N}$ such that $w_n<\epsilon$ on $F_k$, which implies that  $F_k$ is totally bounded, hence compact and the capacity of $(\E,\D(\E))$ is tight. In the following we show \eqref{eq2.8}: we fix $m\in\mathbb{N}$,
consider
$v_m\in  D(\E )$ and
$$D_{h_k}v_m(\gamma)=\varphi'(\bar{d}(\gamma,\xi_m)^2)\sum_{i=1}^N2\int_0^1( \psi^i(\gamma(s))-\psi^i(\xi_m(s)))\left\langle U_s^{-1}\nabla \psi^i(\gamma(s)),h_k(s)\right\rangle \dd s.$$
Thus we obtain
\begin{equation}\label{eq2.9}\aligned\E(v_m,v_m)=&\frac{1}{2}\int_{\bf E}\sum_{k=1}^\infty\bigg(D_{h_k}v_m(\gamma)\bigg)^2\dd\mu\\\leq& 2\int_{\bf E}\sum_{k=1}^\infty \bigg(\sum_{i=1}^N\int_0^1( \psi^i(\gamma(s))-\psi^i(\xi_m(s)))\left\langle U_s^{-1}\nabla \psi^i(\gamma(s)),h_k(s)\right\rangle  ds\bigg)^2\dd\mu\\\leq&C_N\sum_{i=1}^N\int_{\bf E}\|U_\cdot^{-1}\nabla \psi^i(\gamma)\|_\mathbf{H}^2\dd\mu\leq C, \quad \forall m\in \mathbb{N}.\endaligned\end{equation}
Here $C$ is independent of $m$ and in the last inequality we used that $M$ is compact.
Since $\{\xi_m|m\in\mathbb{N}\}$ is dense in ${\bf E}$, $w_n\downarrow0$ on ${\bf E}$ hence in $L^2({\bf E};\mu)$.
By \eqref{eq2.9} and \cite[IV.Lemma 4.1]{MR92} we have $$\E(w_n,w_n)\leq  C,\quad \forall n\in \mathbb{N}.$$
By \cite[I.2.12, III.3.5]{MR92} we obtain that a subsequence of the Cesaro mean of some subsequence of $w_n$ converges to zero $\E$-quasi-uniformly. But since $(w_n)_{n\in\mathbb{N}}$ is decreasing, \eqref{eq2.8} follows.

For any $\gamma\neq \eta\in {\bf E}$ let $\varepsilon:=\bar{d}(\gamma,\eta)>0$. There exists a certain $\xi_N$ such that $\bar{d}(\xi_N,\eta)<\frac{\varepsilon}{4}$ and $\bar{d}(\xi_N,\gamma)>\frac{\varepsilon}{4}$. Let  $v_m(\gamma):=\varphi(\bar{d}(\gamma,\xi_m)^2),m\in \mathbb{N}$ for $\varphi$ as above. Then $\{v_m\}$ separates the points of ${\bf E}$ and (iii) in the definition of quasi-regular Dirichlet froms (cf. \cite{MR92}) follows. Now the results follow immediately.
$\hfill\square$

\subsection{Dirichlet form on loop sapce}
In this subsection,  we construct the quasi-regular Dirichlet form on loop space. To do that, we first need the integration by parts formula with respect to the Brownian bridge measure and this formula does not only depend on  bounds of the Ricci curvature, but also on  the Hessian  of the logarithm heat kernel on $M$.
Fix $o\in M$, the based loop space $L_{o,o}(M)$ over $M$ is defined by
$$L_{o,o}(M):=\left\{\gamma\in C([0,1];M):\gamma(0)=\gamma(1)=o\right\}.$$
Then $L_{o,o}(M)$ is a Polish
space under the uniform distance $d_\infty$.

As in the previous section, we work with the following simple but natural distance on  $L_{o,o}(M)$,
$$\tilde{d}(\gamma,\eta):={\int_0^1\rho(\gamma_s,\eta_s)^2ds}, \quad \gamma,\eta\in L_{o,o}(M).$$   Let ${\bf E}$ be the closure of $L_{o,o}(M)$ in $\{\eta:[0,T]\rightarrow M; \int_0^1\rho(o,\eta_s)^2ds<\infty\}$ with respect to the distance $\tilde{d}$. Then ${\bf E}$ is a Polish space.

Let $\P_{o,o}$ be the Brownian bridge measure on $L_{o,o}(M)$, which can be extended to a Borel measure on ${\bf E}$.
Let $O(M)$ be the orthonormal frame bundle over $M$, and let $\pi: O(M) \rightarrow M$ be the canonical projection.
Let $(\gamma_t)_{\{0\leq t\leq1\}}$ be the coordinate process on $L_{o,o}(M)$, $(\F_t)_{0\leq t\leq1}$ the $\P_{o,o}$-completed natural filtration of $(\gamma_t)$. We set $\F=\F_1$. Then $(\gamma_t)$ is a semimartingale on the stochastic basis $({\bf E},\F,\F_t,\P_{o,o})$. For a given orthonormal frame $u_0\in \pi^{-1}(x)\subset O(M)$, there exists a unique stochastic horizontal lift $(U_t)$ of $(\gamma_t)$, determined by the Levi-Civita connection, such that $U_0=u_0$. Let
\begin{equation}\label{eq2.18}
\begin{cases}
&\dd B_t=U_t^{-1}\circ \dd \gamma_t-U_t^{-1}\nabla \log p_{1-t}(\gamma_t,o)\dd t,\quad t\geq0\\
& B_0=0,
\end{cases}
\end{equation}
where $\circ \dd\gamma_t$ stands for the Stratonovich differential of $\gamma_t$ and $p_t(x,y)$ is the heat kernel of $\frac{1}{2}\Delta$ with $\Delta:=$ Levi-Civita Laplacian on $M$. $(B_t)_{0\leq t\leq 1}$ is an $\mathbb{R}^d$-valued standard Brownian motion.

By Driver's integration by parts formula \cite{Dri94} (see also \cite{H97, CLW11}) we have for $F\in \F C_b^1, h\in \mathbb{H}_0:=\{h\in \H| h(1)=0\},$
\begin{equation}\label{eq2.19}\int_{L_{o,o}(M)} D_hF\dd \P_{o,o}=\int_{L_{o,o}(M)} F\beta_h\dd \P_{o,o},\end{equation}
with
 $$L^2({\bf E},\P_{o,o})\ni\beta_h:=\int^1_0\left\langle {h}'_s+\frac{1}{2}\Ric_{U_s}h_s-\textrm{Hess}_{U_s}\log p_{1-s}(\cdot,o)h_s,\dd B_s\right\rangle,$$
 where $\textrm{Hess}_{u}fa:=u^{-1}\textrm{Hess}f(\pi(u))ua$ for $u\in O(M), a\in\mathbb{R}^d$ and smooth function $f$ on $M$.
 Let $\{h_k\}\subset \mathbb{H}_0$ be  an orthonormal basis in $\HH$ such that $h_k\in\mathbb{H}_0, k\in\mathbb{N}$.
 Similarly as above we easily deduce that the form
$$\E(F,G):=\frac{1}{2}\int_{\bf E} \langle DF, DG\rangle_{\HH}\dd \P_{o,o}=\frac{1}{2}\sum_{k=1}^\infty\int_{\bf E} D_{h_k}F D_{h_k}G \dd \P_{o,o}, \quad F,G\in\F C_b^1$$
is closable and its closure $(\E,\D(\E))$ is a quasi-regular Dirichlet form on $L^2({\bf E};\P_{o,o})=L^2(L_{o,o}(M);\P_{o,o})$.
Moreover, Theorems 2.2-2.4 still hold in this case.

\subsection{Dirichlet form on free path/loop space}
Similar to the above two subsections, in this subsection, we  construct a class of quasi-regular Dirichlet forms on the free path/loop space. Let $\sigma$ be a probability measure on $M$ and $\dd \sigma (x)=v(x)\dd x$ some $C^1$-function $v$ on $M$, and $\P_\sigma$ be the distribution of the Brownian motion/Brownian bridge starting
from $\sigma$ up to time $1$, which is then a probability measure on the free path/loop space:
\begin{equation}\label{eq2.15}\aligned
W(M)=C([0,1]; M)\textrm{ or } L(M)=\cup_{y\in M}L_{y,y}(M).
\endaligned\end{equation}
In fact, we know that
$$\dd \P_\sigma=\int_M\P_y\dd \sigma(y),$$
where $\P_y$ is the law of Brownian motion/Brownian bridge starting at $y$.
Similarly, we define the $L^2$-distance on $W(M)/L(M)$ by
$$\tilde{d}(\gamma,\eta):={\int_0^1\rho(\gamma_s,\eta_s)^2\dd s}, \quad \gamma,\eta\in W(M).$$  Let ${\bf E}$ be the closure of $W(M)/L(M)$ in $\{\eta:[0,T]\rightarrow M; \int_0^1\rho(o,\eta_s)^2ds<\infty\}$ with respect to the distance $\tilde{d}$. Then ${\bf E}$ is a Polish space. $\P_\sigma$ can be extended to a Borel measure on ${\bf E}$.
Choose a sequence of $\{h_k\}\subset \mathbb{H}$ such that it is an orthonormal basis in $\HH$. Then the quadratic form on the free path/loop space is defined by
$$\E(F,G):=\frac{1}{2}\int_{\bf E} \langle DF, DG\rangle_{\HH}\dd\P_\sigma=\frac{1}{2}\sum_{k=1}^\infty\int_{\bf E} D_{h_k}F D_{h_k}G  \dd\P_\sigma, \quad F,G\in\F C_b^1.$$
By the integration by parts formula in \cite{FW05}/\cite[Lemma 4.1]{CLWFL} (and the references therein): for $F\in \F C_b^1, h\in \mathbb{H},$
\begin{equation}\label{eq2.17}\int D_{h}F\dd\P_\sigma=\int F\beta_{h}\dd\P_\sigma,\end{equation}
where
 $$\beta_{h}:=\int^1_0\left\langle {h}_{s}'+\frac{1}{2}\Ric_{U_s}h_{s},\dd B_{s}\right\rangle\textrm{ or }\int^1_0\left\langle {h}'_s+\frac{1}{2}\Ric_{U_s}h_s-\textrm{Hess}_{U_s}\log p_{1-s}(\cdot,o)h_s,\dd B_s\right\rangle,$$
 and $\beta_{h}\in L^2({\bf E},\P_\sigma)$. Here $B$ is the corresponding Brownian motion in $\mathbb{R}^d$. This  implies that
the form
$\E$
is closable, and similarly as above, we can prove that its closure $(\E,\D(\E))$ is a quasi-regular Dirichlet form on $L^2({\bf E};\P_\sigma)=L^2(W(M);\P_\sigma)$.
Moreover, Theorems 2.2-2.4 still hold in this case.

\beg{Remark}\label{r2.6}
Compared to the proof of the closability of the O-U Dirichlet form $\E^{OU}$ on the free path/loop space in \cite{FW05}, our situation is simpler now. This is because the integration by parts formula for O-U Dirichlet form depends on the initial distribution $\sigma$. The present case does not depend on the initial point since now we take the $L_2$-space as the intermediate space.
\end{Remark}

\section{Properties of $L^2$-Dirichlet form on path space}\label{sect3}
In this section, we  study properties of the stochastic heat process $X_t, t\geq0,$  and $L^2$-Dirichlet form $\E$ constructed in Section 2.1. In fact, we establish some functional inequalities associated with $(\E,D(\E))$. As mentioned in Remark 2.1, the results in Section 2 also hold when $M$ is not compact. Therefore, in this section we drop the compactness condition on $M$.

\subsection{Log-Sobolev inequality}
In this subsection, we establish log-Sobolev inequality for the $L^2$-Dirichlet form.

 \beg{thm}\label{T3.1}[Log-Sobolev inequality] Suppose that $\Ric\geq -K$ for $K\in\mathbb{R}$, then the log-Sobolev inequality holds
\begin{equation}\label{eq3.2}\mu(F^2\log F^2)\le 2C(K) \E(F,F),\ \ \ \ F \in \F C^1_{b},
\ \mu(F^2)=1,\end{equation}
where $C(K):=\frac{e^K-1-K}{K^2}\wedge C_0(K)$ with
\begin{equation*}
C_0(K)=\begin{cases}
& \frac{4}{K^2}\left(1-\sqrt{2e^{\frac{K}{2}}-e^{K}}\right),\quad\quad\text{if}~K<0,\\
&\frac{2}{K^2}\left(e^{K}-2e^{\frac{K}{2}}+1\right),\quad\quad\quad\text{if} ~K>0.
\end{cases}
\end{equation*}
\end{thm}

\beg{Remark}\label{r3.2}  \beg{enumerate}
\item[$(i)$]  In fact, Theorem 3.1 has first been proved in \cite{GW06}. Compared to their results, our constant $C(K)$ is smaller. The constant in \cite{GW06}
is given by
\begin{equation*}
\tilde{C}(K)=\begin{cases}
& \frac{4}{K^2}\left(1-\sqrt{2e^{\frac{K}{2}}-e^{K}}\right),\quad\quad\text{if}~2e^{\frac{K}{2}}-e^K>0,\\
&\frac{2}{K^2}\left(e^{K}-2e^{\frac{K}{2}}+1\right),\quad\quad\quad\text{if} ~2e^{\frac{K}{2}}-e^K<0.
\end{cases}
\end{equation*}
Then it is easy to see that $\tilde{C}(K)\geq C_0(K)$ for $K>0$ and $2e^{\frac{K}{2}}-e^K>0$.

Comparing the classic O-U Dirichlet form $\E^{OU}$  and the $L^2$-Dirichlet form $\E$, we note that the log-Sobolev inequality associated to two Dirichlet forms are completely different. The former requires uniform bounds on the Ricci curvature, and the latter only needs  lower bounds of the Ricci curvature.

\item[$(ii)$] According to \cite{W05}, the log-Sobolev inequality implies hypercontractivity of the associated semigroup $P_t$, in particular, the  $L^2$-exponential ergodicity of the process:
$\|P_tf-\int f d\mu\|_{L^2}\leq e^{-t/C(K)}\|F\|_{L^2}.$

\item[$(iii)$] The log-Sobolev inequality  also implies the irreducibility of the Dirichlet form $(\E,\D(\E))$. It's obvious that the Dirichlet form $(\E,\D(\E))$ is recurrent. Combining these two results, by \cite[Theorem 4.7.1]{FOT94}, for any nearly Borel non-exceptional set $B$,
$$\mathbf{P}^z(\sigma_B\circ\theta_n<\infty,\forall n\geq0)=1, \quad \textrm{ for q.e. } z\in {\bf E}.$$
Here $\sigma_B=\inf\{t>0:X_t\in B\}$, $\theta$ is the shift operator for the Markov process $X$, and for the definition of any nearly Borel non-exceptional set we refer to \cite{FOT94}. Moreover by \cite[Theorem 4.7.3]{FOT94} we obtain the following strong law of large numbers: for $f\in L^1(\textbf{E},\mu)$
$$\lim_{t\rightarrow\infty}\frac{1}{t}\int_0^tf(X_s)\dd s=
\int f\dd\mu, \quad \mathbf{P}^{z}-a.s.,$$
for q.e. $z\in E$.

\end{enumerate}
\end{Remark}

\ \newline\emph{\bf Proof of Theorem \ref{T3.1}.}
By \cite{GW06} we have the martingale representation theorem, that is, for $F\in \F C_b^1$,
\begin{equation}\label{eq3.5}F=\EE(F)+\int^1_0\langle H_s^F, \dd B_s\rangle ,\end{equation}
with \begin{equation}\label{eq3.7} H_s^F=\EE\left[M_s^{-1}\int_s^1M_\tau(DF(\tau))\dd\tau\bigg|\F_s\right].\end{equation}
Here and in the following $\EE$ means the expectation w.r.t. $\mu$, $B$ is $\mathbb{R}^d$-valued Brownian motion under $\mu$, $(\F_t)$ is the normal filtration generated by $B$ and $M_t$ is the solution of the equation
\begin{equation}\label{eq3.4}\frac{\dd}{\dd t}M_t+\frac{1}{2}M_t\Ric_{U_t}=0,\quad M_0=I.\end{equation}
Let $F=G^2$ for $G\in \F C_b^1$ and consider the continuous version of the martingale $N_s=\EE[F|\F_s]$. We have
$$N_s=\EE F+\int_0^s\langle H_\tau^F,\dd B_\tau\rangle.$$
Now applying It\^{o}'s formula to $N_s\log N_s$, we have
\begin{equation}\label{eq3.8}\EE N_1\log N_1-\EE N_0\log N_0=\frac{1}{2}\EE\int_0^1N_s^{-1}|H_s^F|^2\dd s.\end{equation}
Here and in the following we use $|\cdot|$ to denote the norm in $\mathbb{R}^d$.
On the other hand,
$$DF=D(G^2)=2GDG.$$
Using this relation in the explicit formula \eqref{eq3.7} for $H^F$, we have
\begin{equation}\label{eq3.9}H_s^F=2\EE\left[GM_s^{-1}\int_s^1M_\tau DG(\tau) \dd\tau\bigg|\F_s\right].\end{equation}
By the lower bound on the Ricci curvature, we have
\begin{equation}\label{eq3.6}\|M_s^{-1}M_\tau\|\leq e^{K(\tau-s)/2}.\end{equation}
By Cauchy-Schwarz inequality in \eqref{eq3.9} and \eqref{eq3.6}, we have
$$|H_s^F|^2\leq 4 \EE[G^2|\F_s]\EE\bigg[\bigg(\int_s^1e^{K(\tau-s)/2}|D_\tau G|\dd\tau\bigg)^2\bigg|\F_s\bigg].$$
Here and in the following we use $D_\tau G$ to denote $DG(\tau)$ for simplicity. Thus the right hand side of \eqref{eq3.8} can be controlled by
\begin{equation}\label{eq3.1}\aligned 2\EE\int_0^1\bigg[\bigg(\int_s^1e^{K(\tau-s)/2}|D_\tau G|d\tau\bigg)^2\bigg]\dd s\leq&2\EE\int_0^1\bigg[\int_s^1e^{K(\tau-s)}d\tau\int_s^1|D_\tau G|^2\dd\tau\bigg]ds
\\\leq& \frac{2}{K}\int_0^1[e^{K(1-s)}-1]\dd s\E(G,G)\\\leq& 2\frac{e^K-1-K}{K^2}\E(G,G).\endaligned\end{equation}
Now we use another way to control the left hand side of \eqref{eq3.1}.  We have the following estimate, which follows essentially from \cite{FW15}: H\"{o}lder's inequality implies that
$$\aligned \bigg(\int_s^1e^{K(\tau-s)/2}|D_\tau G|\dd\tau\bigg)^2\leq \int_s^1e^{K(\tau-s)/2}\dd\tau\int_s^1e^{K(\tau-s)/2}|D_\tau G|^2\dd\tau.\endaligned$$
Then  changing the order of integration we obtain
$$\aligned \EE\int_0^1\bigg[\bigg(\int_s^1e^{K(\tau-s)/2}|D_\tau G|\dd\tau\bigg)^2\bigg]\dd s\leq \EE\int_0^1J_1(s)|D_sG|^2\dd s,\endaligned$$
where
$$\aligned J_1(s)=&-\int_0^s\frac{2}{K}\big(1-e^{K(1-t)/2}\big)e^{K(s-t)/2}\dd t\\=&\frac{2}{K^2}\bigg[2(1-e^{\frac{Ks}{2}})-e^{\frac{K(1-s)}{2}}
+e^{\frac{K(1+s)}{2}}\bigg].\endaligned$$
Taking the derivative of $t\rightarrow J_1(t)$ gives
$$J_1'(t)=-\frac{2}{K}\bigg[e^{\frac{Kt}{2}}-\frac{1}{2}e^{\frac{K(1-t)}{2}}
-\frac{1}{2}e^{\frac{K(1+t)}{2}}\bigg].$$
In addition, we have $J_1(0)=0, J_1(1)=\frac{2}{K^2}\big(1-e^{\frac{K}{2}}\big)^2.$
Moreover, $$J_1'(0)=-\frac{2}{K}\bigg[1-e^{\frac{K}{2}}\bigg]>0, \quad J_1'(1)=\frac{1}{K}\bigg[e^{\frac{K}{2}}-1\bigg]^2.$$
For $J_1'(t)=0$ we only have one solution
$e^{-Kt}=2e^{-\frac{K}{2}}-1$,
which implies that when $K>0$, $J_1'(s)\geq0$ for all $s\geq0$. Then for $K>0$, $J(s)$ is increasing, which implies that  $$\aligned \EE\int_0^1\bigg[\bigg(\int_s^1e^{K(\tau-s)/2}|D_\tau G|\dd\tau\bigg)^2\bigg]\dd s\leq J_1(1)\EE\int_0^1|D_sG|^2\dd s\leq C_0(K)\EE\int_0^1|D_sG|^2\dd s.\endaligned$$
For $K<0$, we suppose $t_0\in(0,1)$ satisfying  $J_1'(t_0)=0$, which is the maximum point of $J_1$.
Then for $K<0$,
$$\aligned \EE\int_0^1\bigg[\bigg(\int_s^1e^{K(\tau-s)/2}|D_\tau G|\dd\tau\bigg)^2\bigg]\dd s\leq J_1(t_0)\int_0^1|D_sG|^2\dd s\leq C_0(K)\int_0^1|D_sG|^2\dd s.\endaligned$$
Combining all the above, we complete the proof.
$\hfill\square$
\vskip.10in

 In the following Theorem \ref{T3.4}, we  obtain a new constant for the log-Sobolev inequality for  Einstein manifolds. In this case the constant $C(K)$ tends to the optimal constant in the flat case as $K\rightarrow0$ (see \cite{CG02}).

\beg{thm}\label{T3.4} Suppose that $M$ is an Einstein manifold with constant Ricci curvature $-K\in \R$. Then the log-Sobolev inequality for
$(\E,\D(\E))$ holds:
\begin{equation}\label{eq3.11}\mu(F^2\log F^2)\le 2C(K) \E(F,F),\ \ \ \ F \in \F C^1_{b},
\ \mu(F^2)=1,\end{equation}
where $$C(K):=\left\{\left[4d\left(\sum_k|A_k|\right)^2\frac{e^K-1}{K}+2A_0^2\right]^{1/2}+\frac{2\pi}{(K^2+\pi^2)}\vee\frac{1}{\pi}\right\}^2$$
and
$$A_k:=\left[\frac{ K}{2}+\frac{2\pi^2}{ K}\left(k+\frac 1 2\right)^2\right]^{-1}, \quad k\in \mathbb{N}\cup\{0\}.$$
\end{thm}

\beg{Remark}
In fact, we have
$$\lim_{K\rightarrow0}C(K)=\frac{4}{\pi^2},$$
and $\frac{4}{\pi^2}$ is the optimal constant in the $\mathbb{R}^d$ case (see \cite{CG02, D04} and the references therein).
\end{Remark}

\begin{proof} Let $ h_{\alpha,k}:= \sqrt{2}\sin\left[
\left(k+\frac 1 2\right)\pi \tau\right] e_\alpha$ for $\alpha=1,...,d, k\in \mathbb{N}\cup\{0\}$. Here $\{e_\alpha\}$ is the usual orthonormal basis for $\mathbb{R}^d$ given by $e_\alpha=(0,...,1,..,0)$. It is easy to see that $\{h_{\alpha,k}\}$ is an orthonormal basis of ${\bf H}$. We start with the following computation:
$$\aligned
\int_s^1e^{K(\tau-s)/2} h_{\alpha,k}\dd\tau
=e_\alpha e^{-Ks/2}\int_s^1e^{K\tau/2} \sqrt{2}\sin\left[
\left(k+\frac 1 2\right)\pi \tau\right] \dd\tau:=\sqrt{2}e_\alpha e^{-Ks/2}B(s,k,K),\endaligned$$
with
$$\aligned B(s,k,K)&=\left[\frac{ K}{2}+\frac{2\pi^2}{ K}\left(k+\frac 1 2\right)^2\right]^{-1}\\
&\times \bigg\{(-1)^ke^{K/2}-e^{Ks/2}\sin\left[
\left(k+\frac 1 2\right)\pi s\right]
+ \frac{2\pi}{ K}\left(k+\frac 1 2\right)e^{Ks/2}\cos\left[
\left(k+\frac 1 2\right)\pi s\right]\bigg\}.\endaligned$$
Thus, we have
\begin{equation}\label{eqA}\aligned
\int_s^1e^{\frac{K(\tau-s)}{2}} h_{\alpha,k}\dd\tau
&=A_k(-1)^ke^{K(1-s)/2}\sqrt{2}e_{\alpha}-A_kh_{\alpha,k}+B_k\bar h_{\alpha,k},
\endaligned\end{equation}
where
$$\aligned&A_k=\left[\frac{ K}{2}+\frac{2\pi^2}{ K}\left(k+\frac 1 2\right)^2\right]^{-1},\quad B_k=\left[\frac{ K}{2}+\frac{2\pi^2}{ K}\left(k+\frac 1 2\right)^2\right]^{-1}\frac{2\pi}{ K}\left(k+\frac 1 2\right),\\&
 \bar h_{\alpha,k}(\tau):=\sqrt{2}\cos\left[
\left(k+\frac 1 2\right)\pi \tau\right] e_\alpha.\endaligned$$
It is easy to see that $B_k\leq \frac{2\pi}{K^2+\pi^2}\vee\frac{1}{\pi}.$  Indeed,  if $K^2<\pi^2$ then $B_k$ is decreasing with respect to $k$ and if $K^2\geq\pi^2$ then $B_k\leq \frac{1}{|K|}\leq \frac{1}{\pi}.$
Now since $\Ric=-K $ we have $M_s^{-1}M_\tau=e^{\frac{K(\tau-s)}{2}}I$. A similar argument as in the proof of Theorem \ref{T3.1} implies the left hand side of \eqref{eq3.1} can be controlled by
$$\aligned
&2\EE\bigg[\int_0^1\bigg|
\sum^d_{\alpha=1}\sum^\infty_{k=0}\left\<DG,h_{\alpha,k}\right\>_{\mathbf{H}} \int_s^1e^{\frac{K(\tau-s)}{2}} h_{\alpha,k}\dd\tau\bigg|^2\dd s\bigg]\\
=&2\EE\bigg[\int_0^1\bigg|
\sum^d_{\alpha=1}\sum^\infty_{k=0}\left\<DG,h_{\alpha,k}\right\>_{\mathbf{H}} A_k(-1)^k\sqrt{2}e^{K(1-s)/2}e_{\alpha}\\
&-\sum^d_{\alpha=1}\sum^\infty_{k=0}\left\<DG,h_{\alpha,k}\right\>_{\mathbf{H}}
A_kh_{\alpha,k}
+\sum^d_{\alpha=1}\sum^\infty_{k=0}\left\<DG,h_{\alpha,k}\right\>_{\mathbf{H}}
B_k\bar h_{\alpha,k}\bigg|^2\dd s\bigg]
\\=&2\bigg(\EE\int_0^1\bigg|
\sum^d_{\alpha=1}\sum^\infty_{k=0}\left\<DG,h_{\alpha,k}\right\>_{\mathbf{H}} A_k(-1)^k\sqrt{2}e^{K(1-s)/2}e_{\alpha}-\sum^d_{\alpha=1}\sum^\infty_{k=0}\left\<DG,h_{\alpha,k}\right\>_{\mathbf{H}}
A_kh_{\alpha,k}\bigg|^2\dd s
\\&+2\EE\int_0^1
\bigg\langle \sum^d_{\alpha=1}\sum^\infty_{k=0}\left\<DG,h_{\alpha,k}\right\>_{\mathbf{H}}A_k(-1)^k\sqrt{2}e^{K(1-s)/2}e_{\alpha}
-\sum^d_{\alpha=1}\sum^\infty_{k=0}\left\<DG,h_{\alpha,k}\right\>_{\mathbf{H}}
A_kh_{\alpha,k},\\&\sum^d_{\alpha=1}\sum^\infty_{k=0}\left\<DG,h_{\alpha,k}\right\>_{\mathbf{H}}
B_k\bar h_{\alpha,k}\bigg\rangle \dd s+\EE\int_0^1\bigg|\sum^d_{\alpha=1}\sum^\infty_{k=0}\left\<DG,h_{\alpha,k}\right\>_{\mathbf{H}}
B_k\bar h_{\alpha,k}\bigg|^2\dd s\bigg)
\\=&:2(I_1+I_2+I_3),\endaligned$$
where we used \eqref{eqA} in the first equality.
Then we have for $I_1$
$$\aligned I_1\leq &4d\EE\int_0^1|DG|^2_{{\bf H}}\bigg(\sum_k|A_k|\bigg)^2e^{K(1-s)}\dd s+2\EE|DG|^2_{{\bf H}}A_0^2
\\=&\EE|DG|^2_{{\bf H}}\bigg[4d(\sum_k|A_k|)^2\frac{e^K-1}{K}+2|A_0|^2\bigg]=:C_1(K)\E(G,G), \endaligned$$
where the first inequality is due to that $A_k^2$ is decreasing w.r.t. $k$ and  $\{h_{\alpha,k}\}$ is an orthonormal basis of ${\bf H}$. For $I_3$ we have
$$\aligned I_3=\sum_{\alpha=1}^d\sum_{k=0}^\infty\left\<DG,h_{\alpha,k}\right\>_{{\bf H}}^2B_k^2
\leq\left(\frac{4\pi^2}{(K^2+\pi^2)^2}\vee\frac{1}{\pi^2}\right)\E(G,G)=:C_2(K)\E(G,G).\endaligned$$
Using H\"{o}lder's inequality we obtain $$\aligned I_2\leq2I_1^{1/2}I_3^{1/2}\leq 2C_1(K)^{1/2}C_2(K)^{1/2}\E(G,G).\endaligned$$
Combining the above estimates we obtain
$$\aligned \mu(F^2\log F^2)
&\leq2(C_1(K)^{1/2}+C_2(K)^{1/2})^2\E(G,G).\endaligned$$
\end{proof}

\subsection{Characterization of the lower bound of the Ricci curvature}

The upper and lower bounds for the Ricci curvature on a Riemannian manifold were well characterized in terms of the twisted Malliavin gradient-Dirichlet form $\E^{OU}$ for the O-U process on the path space (see \textbf{A} in the introduction) in \cite{N,WW16,W1,CT}.
If the Malliavin gradient  is replaced by the $L^2$-gradient $DF$, then we  obtain   characterizations for the lower boundedness of the Ricci curvature in terms of a properly decomposition of the  $L^2$ gradient -Dirichlet form. This subsection is devoted to prove such characterizations.

In fact, all the results in Section 2 and Theorem 3.1 also hold when we change $1$ to  $T>0$. To state our results, let us first introduce some notations: For any point $y\in M$ and $T>0$,
let $x_{y,[0,T]}$ be the Brownian motion starting from $y\in M$ up to $T$, and $\mu_{T,y}$ be the distribution of  $x_{y,[0,T]}$.  Define $\F C_b^T$  as in \eqref{eq2.2} with $1$ replaced by $T$. For any $n\geq1$ and $G\in \F C_b^T$, define the following quadratic form
$$\aligned\E ^K_{T,n,y}(G,G)&=(1+n)C_1(K)\int_{W_y^T(M)}\int_0^{T-\frac{1}{n}}|DG(\gamma)(s)|_{\mathbb{R}^d}^2\dd s\dd  \mu_{T,y}(\gamma)\\
&~~~~~+\left(\frac{1}{n}+\frac{1}{n^2}\right)C_{2,n}(K)\int_{W_y^T(M)}\int_{T-\frac{1}{n}}^T|DG(\gamma)(s)|_{\mathbb{R}^d}^2\dd s\dd\mu_{T,y}(\gamma).\endaligned$$
where
$$C_1(K)=\left[\frac{1}{K^2}\big(TKe^{KT}-e^{KT}+1\big)\right]\bigvee \frac{T^2}{2}, \quad C_{2,n}(K)=\frac{e^{KT}-1}{K}\left(1\vee e^{-\frac{K}{n}}\right).$$
Similarly as in the proof of Theorem \ref{T2.1} we see that $(\E ^K_{T,n,y},\F C_b^T)$ is closable and its closure is a Dirichlet form.
Let $p_t$ be the Markov semigroup of the Brownian motion  $x_y$ starting from $y\in M$, i.e. given by $p_tf(y)=\EE [f(x_{y,t})],y\in M,f\in \B_b(M),t\geq0$.
Let $C_0^\infty(M)$ denote the set of all smooth  functions with compact supports on $M$.
\vskip.10in

 \beg{thm}\label{T3.6}  For  $K\in\mathbb{R}$, the following statements are equivalent:
\beg{enumerate}
\item[$(1)$]
$\Ric\geq -K $.

\item[$(2)$] For any $f\in C_0^\infty(M),T_1>T_2\geq0$ and $y\in M$, we have
$$\left|\int^{T_1}_{T_2}\nabla p_sf(y)\dd s\right| \leq \int^{T_1}_{T_2}\e^{\frac{Ks}{2}} p_s|\nabla f|(y) \dd s .$$

\item[$(3)$]  For any $F\in \left\{\sum_{i=1}^na_i\int_{s_i}^{t_i}f_i(\gamma_s)\dd s,n\in\mathbb{N},f_i\in C_0^\infty(M),s_i,t_i\in [0,T],a_i\in\mathbb{R}\right\}$ and $y\in M, T>0$
$$\left|\nabla_y \EE F(x_{y,[0,T]})\right|\leq \int_0^Te^{\frac{Ks}{2}}\EE|DF(x_{y,[0,T]})|(s)\dd s,$$
where $\nabla_y$ denotes the gradient w.r.t. $y$ and $\EE$ means  expectation w.r.t. $\mu_{T,y}$.

\item[$(4)$]
 For any $y\in M,T>0$, the following log-Sobolev inequality holds for any $n\in\mathbb{N}$:
$$\mu_{T,y}(F^2\log F^2)\leq 2\E^K_{T,n,y}(F,F),\quad F\in\F C^T_b,~\mu_{T,y}(F^2)=1.$$

\item[$(5)$] For any $y\in M,T>0$, the following Poincar\'{e}-inequality holds for any $n\in\mathbb{N}$:
$$\mu_{T,y}(F^2)\leq \E^K_{T,n,y}(F,F),\quad F\in\F C^T_b,~\mu_{T,y}(F)=0.$$

\end{enumerate}
\end{thm}

\proof $(1)\Rightarrow (3)$ By the gradient formula in \cite{H97} (see also \cite{FW05,W14, WW16}), for
$F(\gamma):=f(\gamma_{t})$ with $f\in C_0^\infty(M)$ we have
$$\nabla_y\EE[F(x_{y,[0,T]})]=U_0^y \EE\Big[M_{t}^y(U_{t}^y)^{-1} \nabla f(x_{y,t})\Big],$$
where $U_t^y$ is the solution to \eqref{eq2.1} with $o$ replaced by $y$ and $M_t^y$ is the solution to \eqref{eq3.4} with $U_t$ replaced by $U_t^y$.
Applying the above formula to $F(\gamma)=\sum_{i=1}^na_i\int_{s_i}^{t_i}f_i(\gamma_s)\dd s$, we have
\begin{equation}\label{eq3.10}\aligned\nabla_y\EE[F(x_{y,[0,T]})]&=\sum^n_{i=1}a_i\int^{t_i}_{s_i}\nabla_y\EE[f_i(x_{y,s})]\dd s
\\&=U_0^y\sum^n_{i=1}a_i\int^{t_i}_{s_i}\EE\Big[M_{s}^y(U_{s}^y)^{-1} \nabla f_i(x_{y,s})\Big]\dd s\\
&=U_0^y\int^T_0\EE\big[M_{s}^yDF(x_{y,[0,T]})(s)\big] \dd s
.\endaligned\end{equation}
Combining this with   $\Ric\geq -K$, we have
\begin{equation}\label{eq3.11}|\nabla_y \EE F(x_{y,[0,T]})|\leq \int_0^Te^{\frac{Ks}{2}}\EE|DF(x_{y,[0,T]})|(s)\dd s
.\end{equation}

$(3)\Rightarrow(2)$
Taking $F(\gamma):=\int^{T_1}_{T_2}f(\gamma_s)ds$ for  $f\in C^1_0(M)$ and $T_1>T_2\geq0$, by $(3)$ we have
\begin{equation}\label{eq3.12}\aligned \left|\int^{T_1}_{T_2}\nabla p_sf(y)d s\right|&=\left|\int^{T_1}_{T_2}\nabla_y\EE[f(x_{y,s})]\dd s\right|=\big|\nabla_y\EE[F]\big|\\
&\leq\int^{T}_{0}e^{\frac{Ks}{2}}\EE|DF(x_{y,[0,T]})|(s)\dd s
\\&=\int^{T_1}_{T_2}e^{\frac{Ks}{2}} p_s|\nabla f| (y)\dd s.\endaligned\end{equation}

$(2)\Rightarrow(1)$ Let $f\in C^\infty_0(M)$ with $|\nabla f(y)|=1$ and $\Hess_f(y) = 0$. For any $T>0$ and $\varepsilon>0$, according to $(2)$, we obtain
\begin{equation}\label{eq3.13}
\left|\int^{T+\varepsilon}_T\nabla p_sf(y)d s\right| \leq \int^{T+\varepsilon}_T\e^{\frac{Ks}{2}} p_s|\nabla f|(y) \dd s.
\end{equation}
Dividing the two sides of the above equation  by $\varepsilon$ and letting $\varepsilon$ go to zero, we get
\begin{equation}\label{eq3.14}
\left|\nabla p_Tf\right|(y) \leq \e^{\frac{KT}{2}} p_T|\nabla f|(y).
\end{equation}
Then by the classical result (or refer to \cite{W05} and references therein), (1) follows. Note that $p_{2T}=\tilde p_T$, where  $\tilde p_T$ is the semigroup associated with the generator $\Delta$. Thus we complete the proof of this step.

$(5)\Rightarrow (1)$ Let $f\in C^\infty_0(M)$ with $|\nabla f(y)|=1$ and $\Hess_f(y) = 0$. Taking $F(\gamma)=n\int_{T-1/n}^Tf(\gamma_s)ds,$ then $DF(\gamma)(s)=nU_s^{-1}\nabla f(\gamma_s)1_{s\in[T-1/n,T]}$. By (5) we have
that $$\aligned &\EE\left[n\int_{T-1/n}^Tf(\gamma_s)\dd s\right]^2-\left[\EE n\int_{T-1/n}^Tf(\gamma_s)\dd s\right]^2
\\\leq& C_{2,n}(K)(n+1)\int_{T-1/n}^T\EE|\nabla f(\gamma_s)|^2\dd s.\endaligned$$
Letting $n\rightarrow\infty$ we obtain
\begin{equation}\label{eq PI}p_Tf^2-(p_Tf)^2\leq C_2(K)p_T|\nabla f|^2,\end{equation}
with $C_{2}(K)=\frac{e^{KT}-1}{K}.$ According to \cite[Theorems 3.2.3]{W14},   we have
\beg{equation}\beg{split} \label{RIC}   \Ric(\nn f,\nn f)(y)  &= \lim_{t\downarrow 0} \ff 1 T \bigg(\ff{\tilde p_Tf^2(y)-(\tilde p_Tf)^2(y)}{2T}-|\nn \tilde p_Tf(y)|^2\bigg)\\
&= \lim_{T\downarrow 0} \ff 1 T \bigg(\ff{p_{2T}f^2(y)-(p_{2T}f)^2(y)}{2T}-|\nn p_{2T}f(y)|^2\bigg).\end{split}\end{equation}
Combining the above inequality with \eqref{eq PI}, we get
\beg{equation*}\beg{split}\Ric(\nn f,\nn f)(y)  &\leq \lim_{T\downarrow 0} \ff 1 T \bigg(\ff{e^{2KT}-1}{2KT}p_{2T}|\nabla f|^2(y)-|\nn p_{2T}f(y)|^2\bigg)\\
&=\lim_{T\downarrow 0} \ff 1 T \bigg(\left(\ff{e^{2KT}-1}{2KT}-1\right)p_{2T}|\nabla f|^2(y)\bigg)+\lim_{T\downarrow 0} \ff 1 T \bigg(p_{2T}|\nabla f|^2(y)-|\nn p_{2T}f(y)|^2\bigg)\\
&= K +\lim_{T\downarrow 0} \ff 1 T \bigg(\tilde p_{T}|\nabla f|^2-|\nn \tilde p_{T}f(y)|^2\bigg)
\\
&= K +2\Ric(\nn f,\nn f)(y),\end{split}\end{equation*}
where the last inequality follows due to the formula in \cite[Theorems 3.2.3]{W14}:
$$\Ric(\nn f,\nn f)(y)= \lim_{T\downarrow 0} \ff{\tilde p_T|\nn f|^2(y)-|\nn \tilde p_T f|^2(y)}{2T}.$$
Therefore, we complete the proof of this step.

$(1)\Rightarrow (4)$ According to the proof of Theorem \ref{T3.1}, we only need to prove the following: for any $n\in\mathbb{N}$
$$\EE \int_0^T\left[\int_t^T e^{\frac{(s-t)K}{2}}|DF|\dd s\right]^2\dd t\leq \tilde{\E }_{T,n,y}^K(F,F).$$
In fact, we know that
$$\aligned &\EE \int_0^T\left[\int_t^T e^{\frac{(s-t)K}{2}}|DF|\dd s\right]^2\dd t
\\\leq&\EE \int_0^T\left[\int_{t\vee(T-1/n)}^{T-1/n} e^{\frac{(s-t)K}{2}}|DF|\dd s+\int_{T-1/n}^T e^{\frac{(s-t)K}{2}}|DF|\dd s\right]^2\dd t
\\\leq&(1+n)\int_0^T(T-t)\EE \int_{t\vee(T-1/n)}^{T-1/n} e^{(s-t)K}|DF|^2\dd s \dd t\\
&+(1/n+1/n^2)\int_0^T\EE\int_{T-1/n}^T e^{(s-t)K}|DF|^2\dd s\dd t
\\\leq&(1+n)\int_0^T(T-t)(1\vee e^{(T-t)K})\dd t\EE \int_0^{T-1/n} |DF|^2\dd s
\\&+(1/n+1/n^2)\int_0^Te^{(T-t)K}\dd t(1\vee e^{-K/n})\EE\int_{T-1/n}^T |DF|^2\dd s
\\=&(1+n)C_1(K)\EE \int_0^{T-1/n} |DF|^2\dd s+(1/n+1/n^2)C_2(K,n)\EE\int_{T-1/n}^T |DF|^2\dd s,\endaligned$$
where we used H\"{o}lder's inequality and Young's inequality in the second inequality.
Thus we obtain the result.$\hfill\Box$

\vskip.10in
\section{Stochastic heat equation}
Based on the Andersson-Driver approximation of the Wiener measure, we now present a heuristic derivation of the equation for the process (constructed in Section 2) on path space. When $M$ is Euclidean space, we may choose some suitable linear functions, which are in the domain of the generator, through which we can deduce the associated stochastic heat equation. However, when $M$ is a Riemannian manifold, in general, it is not easy to find suitable test functions on $\mathbf{E}$ belonging to the domain of the generator and  derive the associated equation. Instead,  we will use a suitable approximation to give some intuitive idea how to deduce the equation.  As mentioned in Section 1, it is proved in \cite{AD99} that natural approximations of $\exp(-\frac{1}{2}E(\gamma))\D\gamma$ do indeed converge to Wiener
measure on $M$. For the sake of simplicity, in this section, we suppose that $M$ is compact. First we write the equations associated with the approximation measures.

\subsection{Preliminary}

Before going on, we need to introduce some notations  from \cite{AD99}. We will also use $\langle\cdot,\cdot\rangle$ to denote the Riemmanian metric. Let $\T$ be the set of all partitions of $[0,1]$ and
\begin{equation}\label{eq4.1}E(\gamma):=\int^1_0\langle\gamma'(s),\gamma'(s)\rangle \dd s\end{equation} for all absolutely  continuous curves $\gamma\in W_o(M)$, where $\gamma'(s):=\frac{\dd }{\dd s}\gamma(s)$. Otherwise, set $E(\gamma)=\infty$.
Define the  space of finite energy paths:
\begin{equation}\label{eq4.2}\aligned
H(M):=\left\{\gamma\in W_o(M): \gamma \textrm{ is an absolutely continuous curve and } E(\gamma)<\infty\right\}.
\endaligned\end{equation}
For each $\gamma\in H(M)$,
the tangent space $T_\gamma H(M)$ of $H(M)$ at $\gamma$ may be naturally identified with the space of all absolutely continuous vector fields $X:[0,1]\rightarrow TM$ along $\gamma$ with $X(0)=0$ and $G^1(X,X)<\infty$, where
\begin{equation}\label{eq 4.4}\aligned
G^{1}(X,X):=\int^1_0\left\langle\frac{\nabla X(s)}{\dd s},\frac{\nabla X(s)}{\dd s}\right\rangle\dd s,
\endaligned\end{equation}
\begin{equation}\label{eq 4.5}\aligned
\frac{\nabla X(s)}{\dd s}:=\parals_s(\gamma)\frac{\dd}{\dd s}(\parals_s^{-1}(\gamma)X(s))
\endaligned\end{equation}
and $\parals_s(\gamma):T_oM\rightarrow T_{\gamma(s)}M$ is parallel translation along $\gamma$ relative to the Levi-Civita covariant derivative $\nabla$. As mentioned in \cite{AD99}, on the tangent space $TH(M)$ there exists a natural metric given by
\begin{equation}\label{eq4.4}\aligned
G^0(X,X):=\int^1_0\left\langle X(s),X(s)\right\rangle \dd s,
\endaligned\end{equation}
for any $X\in TH(M)$.

Now we introduce finite dimensional approximations to $(H(M), G^0)$:
for every $\PP:=\{0=s_0<s_1<s_2<...<s_n=1\}\in \T$ with $\Delta_is=s_i-s_{i-1}$, define
\begin{equation}\label{eq4.3}\aligned
H_\PP(M):=\{\gamma\in H(M)\cap C^2([0,1]/ \PP): \nabla \gamma'(s)/ d s=0, s\notin \PP\}.
\endaligned\end{equation}
These are the piecewise geodesics paths in $H(M)$, which change directions only at the partition points. For $\gamma\in H_\PP(M)$ the tangent space $T_\gamma H_\PP(M)$ can be identified with elements $X\in T_\gamma H_\PP(M)$ satisfying the Jacobi equations on $[0,1]\backslash\PP$, see \cite[Prop. 4.4]{AD99} for more details.

By induction, we may easily get the metric on $TH_\PP(M)$ for the partition $\PP\in \T$,
\begin{equation}\label{eq4.5}\aligned
G^0_\PP(X,Y):=\sum^n_{i=1}\left\langle X(s_{i}),Y(s_{i})\right\rangle\Delta_is,
\endaligned\end{equation}
for all $X,Y\in T_\gamma H_\PP(M)$ and $\gamma\in H_\PP(M)$. Let $\textrm{Vol}_{G^0_\PP}$ be the volume form on $H_\PP(M)$ determined by $G^0_\PP$. By the arguments in \cite{AD99},  $\textrm{Vol}_{G^0_\PP}$ may be interpreted as a suitable approximation to $\D\gamma$ mentioned in introduction.

Denote by $\nu^0_\PP$ the measure on $H_\PP(M)$ given by
$$\nu^0_\PP:=\frac{1}{Z_\PP}e^{-\frac{1}{2}E}\textrm{Vol}_{G^0_\PP},$$
where $E:H(M)\rightarrow[0,\infty)$ is the energy functional defined in \eqref{eq4.1} and $Z_\PP$ is a normalization constant given by
$$Z_\PP:=\Pi_{i=1}^n(\sqrt{2\pi}\Delta_is)^d.$$
The following is one of the main results from \cite{AD99}.

\beg{thm}\label{T4.1}  \cite[Theorem 1.8]{AD99}
 Suppose that $f:W_o(M)\rightarrow\mathbb{R}$ is  bounded and continuous,
$$\lim_{|\PP|\rightarrow0}\int_{H_\PP(M)}f(\gamma)\dd\nu^0_\PP(\gamma)=\int_{W_o(M)}f(\gamma)e^{-\frac{1}{6}\int_0^1\textrm{\Scal}(\gamma(s))\dd s}\dd\mu(\gamma),$$
where $\Scal$ is the scalar curvature of $M$ and $\mu$ is the law of Brownian motion on $M$ introduced in Section 2.1.
\end{thm}

For technical reasons, we need to introduce the following  subspace $H_\PP^\delta(M)$ of $H_\PP(M)$ such that every element $\gamma\in H_\PP^\delta(M)$  is a piecewise geodesic and each part $\gamma([s_{i-1},s_i])$ is the unique geodesic linking $\gamma(s_{i-1})$ and $\gamma(s_{i})$ (see \cite[Sec. 5]{AD99}).
In fact, for any partition
$\T \ni\PP:=\{0=s_0<s_1<s_2<...<s_n=1\}$ with $\Delta s_i=\varepsilon$ for $i=1,..,n$ and each
$\delta>0$ less than the injectivity radius of $M$, define
$$H_\PP^\delta(M):=\left\{\gamma\in H_\PP(M):\int_{s_{i-1}}^{s_i}|\gamma'(s)|\dd s<\delta \textrm{ for } i=1,2,...,n\right\},$$
where $s_0=0$. In the following we always suppose that $\delta>0$ is less than the injectivity radius of $M$.
Then we can easily check that $H_\PP^\delta(M)$ is a locally compact separable metric space with the distance given by $$d^\PP(\gamma,\eta)
=\varepsilon\sum_{i=1}^n\rho(\gamma_{s_i},\eta_{s_i}),\quad \forall~\gamma,\eta\in H_\PP^\delta(M).$$
Moreover, it is easy to show that each $\gamma\in H_\PP^\delta(M)$ is determined uniquely by finite points $o, \gamma(s_1),\gamma(s_2),...,\gamma(s_n)$ (see e.g. \cite[Section 5]{AD99}).
By Theorem \ref{T4.1} and \cite[(6.1), Prop. 5.13]{AD99} we have the following convergence result.

\beg{thm}\label{T4.2} 
Suppose that $f:W(M)\rightarrow\mathbb{R}$ is  bounded and continuous.
Then
$$\lim_{|\PP|\rightarrow0}\int_{H_\PP^\delta(M)}f(\gamma)\dd\nu^0_\PP(\gamma)=\int_{W_o(M)}f(\gamma)e^{-\frac{1}{6}\int_0^1\textrm{\Scal}(\gamma(s))\dd s}\dd\mu(\gamma),$$
where $\Scal$ is the scalar curvature of $M$ and $\mu$ is the law of Brownian motion on $M$ introduced in Section 2.1.
\end{thm}

Next, we will recall some basic geometrical concepts of a Riemannian manifold $M$. As in Section 2, let $O(M)$ be the orthonormal frame bundle over $M$ and let $\pi:O(M)\rightarrow M$ denote the bundle of orthogonal frames on $M$.
Let $\X(M)$ be the set of all smooth vector fields and let $\nabla$ be the Riemannian connection on $M$. The curvature tensor is given in terms of the Riemannian connection
$\nabla$  by the following formula:
$$R(X,Y)Z=\nabla_X\nabla_YZ-\nabla_X\nabla_YZ-\nabla_{[X,Y]}Z,$$
for any vector fields $X, Y, Z\in \X(M)$ on $M$, where $ [X,Y]$ is the Lie bracket of vector fields $X$ and $Y$.

The Ricci curvature may be interpreted as the trace of the curvature tensor and the scalar curvature may be considered as the trace of the Ricci curvature tensor on $M$, that is to say,
$$\aligned &\Ric(X):=\sum_{i=1}^d R(X,\bar{e}_i)\bar{e}_i,\quad \forall~X\in T(M),\\
&\Scal(x):=\sum_{i=1}^d \langle\Ric(\bar{e}_i),\bar{e}_i\rangle_{T_xM}, x\in M,\endaligned$$
where $\{\bar{e}_i\}$ is an orthonormal frame. Denote the curvature form by
 $$\Omega(\eta_1,\eta_2)=u^{-1}R(\pi_*\eta_1,\pi_*\eta_2)u,$$
 for all $u\in O(M)$ and $\eta_1,\eta_2\in T_uO(M)$,  and for $a,b\in \mathbb{R}^d$, let
 $$\Omega_u(a,b):=u^{-1}R(ua,ub)u.$$ Define
 $$R_u(v,w)=u^{-1}R(v,w)u,\quad u\in O(M), ~v,w\in T_{\pi(u)}M.$$
Fix $\gamma\in H(M)$ and $X\in T_\gamma H(M)$, define $q_s(X)$ by
 \begin{equation}\label{eq4.6}q_s(X)=\int_0^sR_{u(r)}(\gamma'(r),X(r))\dd r,\end{equation}
 where $u=\paral(\gamma)$ is the horizontal lift of $\gamma$.

The development map $\phi:H(\mathbb{R}^d)\rightarrow H(M)$ is defined by $\phi(b):=\gamma\in H(M)$ for $b\in H(\mathbb{R}^d)$, where $\gamma$ solves the functional differential equation,
$$\gamma'(s)=\paral_s(\gamma)b'(s),\quad \gamma(0)=o.$$
The anti-development map $\phi^{-1}:H(M)\rightarrow H(\mathbb{R}^d)$ is given by $b=\phi^{-1}(\gamma)$,
where $$b(s)=\int_0^s\paral_r^{-1}(\gamma)\gamma'(r)\dd r.$$
For each $h\in C^\infty(H(M)\rightarrow H(\mathbb{R}^d))$ and $\gamma\in H(M)$, let $X^h(\gamma)\in T_\gamma H(M)$ be given by
$$X_s^h(\gamma):=\paral_s(\gamma)h_s(\gamma)\quad \textrm{ for all } s\in [0,1],$$
where $h_s(\gamma):=h(\gamma)(s)$.
Given $\gamma\in H_\PP(M)$, let $H_{\PP,\gamma}$ be the subspace of $H(\mathbb{R}^d)$ given by
$$H_{\PP,\gamma}:=\{v\in H(\mathbb{R}^d): v''(s)=\Omega_{u(s)}(b'(s),v(s))b'(s), \forall s\notin \PP\},$$
where $u=\paral(\gamma)$ and $b=\phi^{-1}(\gamma)$. By \cite{AD99} we know that
$v\in H_{\PP,\gamma}$ if and only if $X^v(\gamma):=\paral(\gamma)v\in T_\gamma H_\PP(M)$.

\subsection{The approximation Dirichlet form $\E^\PP$}

In this subsection we will mainly derive the Dirichlet form associated with the approximation measures $\nu^0_\PP$. To do that, we need to construct a family of special basis on $TH_\PP(M)$.

For any $\varepsilon>0$, take $\T\ni\mathscr{P}=\{0=s_0<s_1<...<s_n=1\}$ with $\Delta s_i=s_i-s_{i-1}=\varepsilon$ for $i=1,...,n$. Let $\{e_a\}$ be an orthonormal basis for $\mathbb{R}^d$ given by $e_a=(0,..,1,...0)$.
Consider the space $l^2(\mathscr{P};\mathbb{R}^d):=\left\{\hat h:\PP\mapsto \R^d\Big|  \|\hat{h}\|^2_{l^2(\mathscr{P};\mathbb{R}^d)}<\infty\right\}$ under the norm given by
$$\left\langle\hat h_1, \hat h_2\right\rangle_{l^2(\mathscr{P};\mathbb{R}^d)}=\varepsilon \sum_{i=1}^n \langle \hat{h}_1(s_i),\hat{h}_2(s_i)\rangle_{{\mathbb{R}^d}}.$$
Choose an orthonormal basis $\hat{h}_{a,i}\in l^2(\mathscr{P};\mathbb{R}^d),  i=1,...,n, a=1,...,d$,  be   given by
$$\hat{h}_{a,i}(s_j)=\begin{cases}&0, ~~\quad\quad\quad\quad\quad j\neq i\\
&\frac{1}{\sqrt{\varepsilon}}e_a,\quad\quad\quad\quad j= i.\end{cases}$$

For fixed $i=1,...,n, a=1,...,d$, define $h_{a,i}:H_\PP(M)\rightarrow H(\mathbb{R}^d)$ by requiring $h_{a,i}(\gamma)\in H_{\PP,\gamma}$ for all $\gamma\in H_{\PP}(M)$ and for $s\in \PP$,
${h}_{a,i}(\gamma)(s)=\hat{h}_{a,i}(s)$ for all $\gamma\in H_{\PP}(M)$.
For $\gamma\in H_{\PP}^\delta(M)$, $h_{a,i}(\gamma)$ is uniquely determined by the above properties (see the proof of Lemma A.1 below). The following lemma  is used to prove the quasi-regularity of the
 approximation Dirichlet form $\E^\PP$.

\beg{lem}\label{L4.3} $\sup_{r\in[0,1]}|h_{a,i}(\gamma)(r)|\in  L^p(H_\PP^\delta(M),\nu^0_\PP), p>1,$ with $\delta>0$ satisfying
 $\cosh(\sqrt{\kappa_0}\delta)\kappa_0\delta^2<1$ .  Here $\kappa_0$ is an upper bound for the norms of the curvature tensor $R$ (or equivalently $\Omega$).
\end{lem}

\begin{proof} We only consider $h_{a,1}(r)$ on $[0,\varepsilon]$. The other cases can be handled similarly.
We use the following notations: $\gamma\in H_\PP^\delta(M)$, $b:=\phi^{-1}(\gamma)$ the anti-development map, $u:=\paral(\gamma)$  and $A(s):=\Omega_{u(s)}(b'(s),\cdot)b'(s)$, $b'(s)=\Delta_i b/\varepsilon$ for $s\in (s_{i-1},s_i]$ with $\Delta_i b=b(s_i)-b(s_{i-1})$.
A similar argument as in the proof of Lemma 8.2 in \cite{AD99} implies that for $r\in [0,\varepsilon]$,
\begin{equation}\label{eq4.7}|h_{a,1}(r)|\leq |h'_{a,1}(0)|\frac{\sinh \sqrt{K}r}{\sqrt{K}},\end{equation}
where
\begin{equation}\label{eq 4.7}K:=\sup_{s\in[0,\varepsilon]}\|A(s)\|\leq \kappa_0\frac{|\Delta_1b|^2}{\varepsilon^2}\leq \kappa_0\frac{\delta^2}{\varepsilon^2}\end{equation}
 and $\|\cdot\|$ is the norm of the matrix. In fact, by Taylor's theorem we have for $s\in [0,\varepsilon]$
\begin{equation}\label{eq4.8}\aligned h(s)=&h(0)+sh'(0)+\int_0^sh''(u)(s-u)\dd u
\\=&sh'(0)+\int_0^sA(u)h(u)(s-u)\dd u.\endaligned\end{equation}
Here and in the following we omit the subindex of $h$ if there's no confusion. Then for $s\in [0,\varepsilon]$
$$|h(s)|\leq s|h'(0)|+K\int_0^s|h(u)|(s-u)\dd u=:f(s).$$
Note that $f(0)=0,f'(s)=|h'(0)|+K\int_0^s|h(u)|\dd u$ and
$$f''(s)=K|h(s)|\leq Kf(s),$$
that is $$f''(s)=Kf(s)+\eta(s),\quad f(0)=0, \quad f'(0)=|h'(0)|,$$
where $\eta(s):=f''(s)-Kf(s)\leq 0.$ Then by the variation of parameter ( cf. the proof of \cite[Lemma 8.2]{AD99}) we have
$$\aligned f(s)=&|h'(0)|\frac{\sinh\sqrt{K}s}{\sqrt{K}}+\int_0^s\frac{\sinh\sqrt{K}(s-r)}{\sqrt{K}}\eta(r)\dd r
\leq |h'(0)|\frac{\sinh\sqrt{K}s}{\sqrt{K}}, \endaligned$$which implies \eqref{eq4.7}.
Also \eqref{eq4.8} implies that
$$h'(0)=\frac{1}{\varepsilon}\left[h(\varepsilon)-\int_0^\varepsilon(\varepsilon-u)A(u)h(u)\dd u\right].$$
Then by \eqref{eq4.7} we have
$$\aligned |h(s)|\leq &\frac{\sinh\sqrt{K}s}{\sqrt{K}}\frac{1}{\varepsilon}\left[\frac{1}{\sqrt{\varepsilon}}+K\int_0^\varepsilon(\varepsilon-u)|h(u)|\dd u\right]
\\\leq&\frac{\sinh\sqrt{K}\varepsilon}{\sqrt{K}}\frac{1}{\varepsilon\sqrt{\varepsilon}}+\frac{\sqrt{K}\varepsilon}{2}\sinh\sqrt{K}\varepsilon\sup_{u\in[0,\varepsilon]}|h(u)|\endaligned $$
By \eqref{eq 4.7} we have $$\sqrt{K}\varepsilon\sinh\sqrt{K}\varepsilon\leq \sqrt{\kappa_0}|\Delta_1b|\sinh(\sqrt{\kappa_0}|\Delta_1b|)\leq {\kappa_0}\delta^2\cosh(\sqrt{\kappa_0}\delta)\leq1.$$
Thus we know that $$\sup_{s\in[0,\varepsilon]}|h(s)|\leq \frac{2}{\varepsilon\sqrt{\varepsilon}}\frac{\sinh\sqrt{K}\varepsilon}{\sqrt{K}}\leq \frac{2}{\sqrt{\varepsilon}}\cosh\sqrt{K}\varepsilon\leq \frac{2}{\sqrt{\varepsilon}}\cosh\kappa_0\delta,$$
which implies the result. Here we used the elementary inequality $\sinh(a)/a\leq \cosh(a)$.
\end{proof}

In the following we fix a $\delta$ as in Lemma \ref{L4.3} and we consider $H_\PP^\delta(M)$ as the state space for the approximation Dirichlet form.  Let
$$\aligned \F C^{\mathscr{P}}_0:=\left\{H_{\mathscr{P}}^\delta(M)\ni\gamma\mapsto F(\gamma):=f(\gamma_{s_1},\gamma_{s_2},...,\gamma_{s_n}),
f\in C_b^1(M^n)\right\}\cap C_0^1(H_\PP^\delta(M)),\endaligned$$
with $C_0^1(H_\PP^\delta(M))$ being continuous, differentiable functions from $H_\PP^\delta(M)$ to $\mathbb{R}$ with compact support.
Since $\gamma\in H_\PP^\delta(M)$, $\gamma$ is determined by $\gamma(s_1),...,\gamma(s_n)$. This implies that every $u=f(\gamma_{t_1},...,\gamma_{t_m})$ with $f\in C_b^1(M^m), 0<t_1<t_2<...<t_m\leq 1$, can be expressed as $g(\gamma_{s_1},\gamma_{s_2},...,\gamma_{s_n})\in \F C^{\mathscr{P}}_0$. By this we can easily conclude that $\F C^{\mathscr{P}}_0$ is dense in $L^2(H_\PP^\delta(M),\nu_\PP^0)$.

For each $F\in \F C^{\mathscr{P}}_0$, the directional derivative of $F$ with respect to $h_{a,i}$ is given by
\begin{equation}\label{eq4.9}D_{h_{a,i}}F(\gamma)=
 \left\langle \nabla_i f(\gamma),\paral_{s_i}(\gamma) h_{a,i}(\gamma)(s_i)
 \right\rangle_{T_{\gamma_{s_i}}M},\quad \gamma\in H_\PP^\delta(M),\end{equation}
where $\nabla_if(\gamma)=\nabla_if(\gamma_{s_1},...,\gamma_{s_n})$. Define for $\gamma\in H_\PP^\delta(M)$,
$$DF(\gamma):=\sum_{i=1}^n\sum_{a=1}^d D_{h_{a,i}}F(\gamma) \hat{h}_{a,i}\in l^2(\PP;\mathbb{R}^d).$$
In this section we also use the notation $DF$ as in Section 2 for simplicity.

\beg{Remark}
By the definition of $h$ and $\hat{h}$ we know that
$DF(\gamma)=\sum_{i=1}^n\sum_{a=1}^d D_{h_{a,i}}F(\gamma)$ ${h}_{a,i}|_{\PP}.$
 For $F\in \F C^{\mathscr{P}}_0$ the directional derivative should be along ${h}_{a,i}\in H_{\PP,\gamma}$, which does not form a basis for $L^2([0,1];\mathbb{R}^d)$. Therefore, we replace $L^2([0,1];\mathbb{R}^d)$ in $\E$ by $l^2(\PP;\mathbb{R}^d)$ in the Dirichlet form $\E^{\mathscr{P}}$ and consider $DF(\gamma)\in l^2(\PP;\mathbb{R}^d)$. For $F\in \F C_b^1, \gamma\in H_\PP^\delta, $ we can find $F_\varepsilon\in \F C^{\mathscr{P}}_0$ such that $\|DF_\varepsilon(\gamma)\|_{l^2(\PP,\mathbb{R}^d)}\rightarrow \|DF(\gamma)\|_{L^2([0,1];\mathbb{R}^d)}$, as $ \varepsilon\rightarrow0$, where the second $DF$ is the $L^2$-gradient in Section 2.
\end{Remark}

Next, we will introduce the quadratic form on $H_{\mathscr{P}}^\delta(M)$.
For any $u,v\in\F C^{\mathscr{P}}_0$, define
$$\E^{\mathscr{P}}(u,v)=\frac{1}{2}\int_{H_{\mathscr{P}}^\delta(M)} \langle Du,Dv\rangle_{l^2(\PP,\mathbb{R}^d)} \dd\nu^0_\PP=\sum_{i=1}^n\sum_{a=1}^d \frac{1}{2}\int_{H_{\mathscr{P}}^\delta(M)} D_{h_{a,i}}u(\gamma)D_{h_{a,i}}v(\gamma)\dd\nu^0_\PP.$$

To prove the closability of the form $\E^\PP$, we need to establish the following integration by parts formula for $\nu^0_\PP$.

\beg{lem}[Integration by parts formula]  For every  $h_{a,j}, a=1,2,...,d, j=1,2,...,n$, we have the following integration by parts formula
\begin{equation}\label{eq4.10}\int_{H_{\mathscr{P}}^\delta(M)}X^{h_{a,j}}f\dd\nu^0_\PP=\int_{H_{\mathscr{P}}^\delta(M)} f\beta_\PP(h_{a,j})\dd\nu^0_\PP\end{equation}
for all $f\in C_0^1(H_{\mathscr{P}}^\delta(M))$ with $\delta$ as in Lemma \ref{L4.3}, where for $p>1$
$$\aligned L^p(H_\PP^\delta(M),\nu_\PP^0)\ni\beta_\PP(h_{a,j})&=\frac{1}{\varepsilon}\left\langle \Delta_jb-\Delta_{j+1}b, h_{a,j}(s_{j}) \right\rangle +\varepsilon\sum_{a_1=1}^d\left\langle
q(X^{h_{a_1,j}})h_{a_1,j},h_{a,j}\right\rangle(s_{j}).
\endaligned$$
Here $b=\phi^{-1}(\gamma)$ with $\Delta_j b=b(s_j)-b(s_{j-1})$ for $j=1,...,n$, $\Delta_{n+1} b=0$ and $q$ is defined in \eqref{eq4.6}.
\end{lem}

\beg{proof} 
 By Stoke's theorem we have for $f\in  C_0^1(H_{\mathscr{P}}^\delta(M))$
 $$0=\int_{H_{\mathscr{P}}^\delta(M)} \left[(X^{h_{a,j}}f)\nu_{\mathscr{P}}^0+fL_{X^{h_{a,j}}}\nu_{\mathscr{P}}^0\right],$$
where we recall
$$\nu_{\mathscr{P}}^0=\frac{1}{Z_\PP^0}e^{-\frac{1}{2}E} \textrm{Vol}_{G^0_\PP}.$$
By the same arguments as in \cite[Lemma 7.3]{AD99} and \eqref{eq4.5}, we know that $\{X^{h_{a,i}},i=1,...,n,a=1,...,d\}$ is a globally defined orthonormal frame for $(H_\PP(M),G^0_\PP)$.  Then
we have for $a=1,...,d, j=1,...,n$,
$$\aligned  L_{X^{h_{a,j}}}\nu_{\mathscr{P}}^0&=-\frac{1}{2} (X^{h_{a,j}}E)(\gamma)\cdot \nu_{\mathscr{P}}^0
+L_{X^{h_{a,j}}} \textrm{Vol}_{G^0_\PP}\\
&=-\frac{1}{2} (X^{h_{a,j}}E)(\gamma)\cdot \nu_{\mathscr{P}}^0
+\frac{1}{Z_\PP^0}e^{-\frac{1}{2}E} \sum_{i=1}^n\sum_{a_1=1}^d G^0_\PP([X^{h_{a_1,i}}, X^{h_{a,j}}], X^{h_{a_1,i}})\cdot \textrm{Vol}_{G^0_\PP}\\
&=-\frac{1}{2} (X^{h_{a,j}}E)(\gamma)\cdot \nu_{\mathscr{P}}^0
+\sum_{i=1}^n\sum_{a_1=1}^d G^0_\PP([X^{h_{a_1,i}}, X^{h_{a,j}}], X^{h_{a_1,i}})\cdot\nu_{\mathscr{P}}^0.\endaligned$$
By the Cartan development map and \cite[Lemma 7.1]{AD99}, we know
$$\aligned(X^{h_{a,j}}E)(\gamma)&=2\int^1_0\left\langle\gamma'(s), \frac{\nabla X^{h_{a,j}}(\gamma)(s)}{\dd s} \right\rangle_{T_{\gamma_{s}}M}\dd s =2\int^1_0\left\langle \paral_s(\gamma)b'(s), \paral_s(\gamma) h_{a,j}'(s) \right\rangle_{T_{\gamma_{s}}M}\dd s\\
&=2\int^1_0\left\langle b'(s), h_{a,j}'(s) \right\rangle\dd s=\frac{2}{\varepsilon}\left\langle \Delta_jb-\Delta_{j+1}b, h_{a,j}(s_{j}) \right\rangle ,\endaligned$$
where we used $b'(s_i)=b'(r)=\Delta_ib/\varepsilon$  for $r\in (s_{i-1},s_i], i=1,...,n,$ in the last equality.

Furthermore, by \cite[Theorem 3.5]{AD99} we have
$$\aligned
&\sum_{i=1}^n\sum_{a_1=1}^dG^0_\PP([X^{h_{a_1,i}}, X^{h_{a,j}}], X^{h_{a_1,i}})\\
=&\varepsilon\sum_{i=1}^n\sum_{a_1=1}^d\sum_{k=1}^n
\left\langle X^{h_{a_1,i}}h_{a,j}-X^{h_{a,j}}h_{a_1,i} ,h_{a_1,i}\right\rangle(s_{k})\\
&-\varepsilon\sum_{i=1}^n\sum_{a_1=1}^d\sum_{k=1}^n\left\langle
q(X^{h_{a,j}})h_{a_1,i}-q(X^{h_{a_1,i}})h_{a,j}, h_{a_1,i}\right\rangle(s_{k})\\
=&\varepsilon\sum_{a_1=1}^d\left\langle
q(X^{h_{a_1,j}})h_{a,j},h_{a_1,j}\right\rangle(s_{j}).
 \endaligned$$
Here we used $X^{h_{a_1,i}}h_{a,j}(s_k)=0$, since $h_{a,j}(s_k)$ is independent of $\gamma$ and we also used $\left\langle
q(X^{h_{a_1,i}})h_{a,j}, h_{a_1,i}\right\rangle(s_{k})\neq0$ only for $i=j=k$ and the skew symmetry of $q(X^{h_{a_1,j}})$ to deduce $\left\langle
q(X^{h_{a,j}})h_{a_1,i},h_{a_1,i}\right\rangle=0$. Thus,
by Stoke's theorem we know that \eqref{eq4.10} holds.

\end{proof}

Based on the above integration by parts formula,  we  obtain the closablity of the following  quadratic form $(\E^\PP, \F C^{\mathscr{P}}_0)$ on $L^2(H_\PP^\delta(M);\nu^0_\PP)$.   Now we prove:

\beg{thm}  The quadratic form $(\E^\PP, \F C^{\mathscr{P}}_0)$
is closable in $L^2(H_\PP^\delta(M);\nu^0_\PP)$ and its closure $(\E^\PP,\D(\E^\PP))$ is a quasi-regular Dirichlet form.
\end{thm}

\beg{proof} {\bf$(a)$ Dirichlet form:}
First we prove that $(\E^\PP,\F C^{\mathscr{P}}_0)$ is closable.
Let $\{F_k\}_{k=1}^{\infty}\subseteq \F C_0^\PP$ with
\begin{equation}\label{eq4.11}
\begin{split}
\lim_{k \rightarrow \infty}\nu^0_\PP\left[F_k^2\right]=0,\ \
\lim_{k,m \rightarrow \infty}\E^\PP\left(F_k-F_m,F_k-F_m\right)=0.
\end{split}
\end{equation}
By \eqref{eq4.11}, we know that $\{{D F_k}\}_{k=1}^{\infty}$ is a Cauchy sequence in
$L^2\left(H_\PP^\delta(M)\rightarrow {l^2(\mathscr{P};\mathbb{R}^d)};\nu^0_\PP\right)$, for which there exists a limit $\Phi$. It suffices to  prove that $\Phi=0$.
Taking $F\in \F C^{\mathscr{P}}_0$, we have for $a=1,...,d, i=1,...,n$,
\begin{equation}\label{eq4.12}\aligned
X^{h_{a,i}} F&=D_{h_{a,i}}F=\langle{D F}, \hat{h}_{a,i}\rangle_{l^2(\mathscr{P};\mathbb{R}^d)}.\endaligned\end{equation}
Thus,
by the above integration by parts formula \eqref{eq4.10}, we have for $G\in \F C_0^\PP$
\begin{equation}\label{eq4.13}
\begin{split}
&\nu^0_\PP\left[G\langle {D F_k}, \hat{h}_{a,i}\rangle_{{l^2(\mathscr{P};\mathbb{R}^d)}}\right]=\nu^0_\PP\left[GX^{h_{a,i}} F_k\right]\\
&=\nu^0_\PP\left[F_kG \beta_\PP(h_{a,i})\right]-\nu^0_\PP\left[F_kX^{h_{a,i}}G\right].
\end{split}
\end{equation}
Since $G$ and $D G$ are bounded and $\beta_\PP(h_{a,i})\in L^2(H_\PP^\delta(M);\nu^0_\PP)$, $F_k$ converges to $0$ in $L^2(\nu^0_\PP)$. By the dominated convergence theorem, taking the limit in \eqref{eq4.13}, we obtain
\begin{equation*}
\begin{split}
&\nu^0_\PP\left[G\langle \Phi, \hat{h}_{a,i}\rangle_{{l^2(\mathscr{P};\mathbb{R}^d)}}\right]=0,\quad \forall\
G \in \F C_0^\PP.
\end{split}
\end{equation*}
Therefore, there exists a $\nu^0_\PP$-null set $\Omega_0 $, such that
$$
\langle \Phi(\gamma), \hat{h}_{a,i}\rangle_{{l^2(\mathscr{P};\mathbb{R}^d)}}=0,\gamma \notin \Omega_0.
$$
Since $\{\hat{h}_{a,i}\}$ is an orthonormal basis in ${l^2(\mathscr{P};\mathbb{R}^d)}$, we  conclude that $\Phi=0$, a.s., and hence
$(\E^\PP,\F C^{\mathscr{P}}_0)$ is closable. Moreover, it is standard that the closure $(\E^\PP,D(\E^\PP))$ is a Dirichlet form.

{\bf$(b)$ Quasi-regularity:} Since  $\gamma\in H_\PP^\delta(M)$ is uniquely determined by
$(\gamma(s_1),\gamma(s_2),...,$ $\gamma(s_n))$, we can easily find a countable dense subset in $\F C_0^\PP$ to separate the points in $H_\PP^\delta(M)$. In fact, similarly as in the proof of Theorem \ref{T2.1}, we use $\psi$ to denote the Nash embedding map. For $k\in\mathbb{N}$, choose $\chi_k\in \F C_0^\PP$ satisfying $\chi_k(\gamma)=1$ if $d(\gamma_{s_{i-1}},\gamma_{s_i})\leq \delta-\frac{1}{k}$ for every $i=1,...,n$. Since $H_\PP^\delta(M)$ is separable we can choose a fixed countable dense set $\{\xi^m|m\in\mathbb{N}\}$ in  $H_\PP^\delta(M)$. Take $\{v_{mki}(\gamma)=|\psi(\gamma_{s_i})-\psi(\xi^m_{s_i})|^2\chi_k(\gamma), k,m\in\mathbb{N}, i=1,...,n\}$, which is a countable dense subset in $\F C_0^\PP$ and separate the points in $H_\PP^\delta(M)$.
Since $H_\PP^\delta(M)$ is locally compact, the tightness of the corresponding capacity follows immediately. Now the quasi-regularity of the Dirichlet form follows.
\end{proof}

Similarly as in Section 2, we can construct a Markov process associated with the above Dirichlet form. We consider $H_{\PP,\Delta}^\delta(M)$ as the one point compactification of $H_\PP^\delta(M)$ (c.f. \cite[P88]{MR92}). Any function $f: H_\PP^\delta(M)\rightarrow\mathbb{R}$ is considered as a function on $H_{\PP,\Delta}^\delta(M)$ by setting $f(\Delta)=0$. By the above proof for quasi-regularity and \cite[Chap. V Corollary 2.16]{MR92} we obtain:

\beg{thm}\label{T4.6}  There exists a  Markov (Hunt) diffusion process
$M^\PP=\big(\Omega,\F,\M_t,$ $(x^\PP_t)_{t\geq0},$ $({P}^z)_{z\in H_{\PP,\Delta}^\delta(M)}\big)$ with state space $H_\PP^\delta(M)$ \emph{properly associated with} $(\E^{\PP},\D(\E^{\PP}))$, i.e. for $u\in L^2(H_\PP^\delta(M);\nu^0_\PP)\cap\B_b(H_\PP^\delta(M))$, the transition semigroup $P_t^\PP u(z):={E}^z[u(x^\PP
 _t)]$ is an $\E^{\PP}$-quasi-continuous version of $T_t^\PP u$ for all $t >0$, where $T_t^\PP$ is the semigroup associated with $(\E^{\PP},\D(\E^{\PP}))$.
\end{thm}
By the integration by parts formula in Lemma 4.5 we can write the explicit martingale solution to the Markov process constructed for $\nu^0_\PP$.

\beg{thm}\label{T4.7}  There exists a \emph{properly  $\E^{\PP}$-exceptional set} $S\subset E$, i.e. $\nu^0_\PP(S)=0$ and ${P}^z[x^\PP(t)\in H_{\PP,\Delta}^\delta(M)\setminus S, \forall t\geq0]=1$ for $z\in H_{\PP,\Delta}^\delta(M)\backslash S$, such that $\forall z\in H_\PP^\delta(M)\backslash S$ under ${P}^z$,  the sample paths of the associated  process $M^\PP$  satisfy the following for  $u(\gamma)=f(\gamma_{s_1},..,\gamma_{s_n})\in \F C_0^\PP$ with $f\in C^\infty(M^n)$,
\begin{equation}\label{eq4.14}
\begin{split}
&u(x_{t}^\PP)-u(x_{0}^\PP)
=\frac{1}{2}\sum_{i=1}^n\sum_{a=1}^d\int_0^tX^{{h_{a,i}}}X^{{h_{a,i}}}u(x_l^\PP)\dd l
\\&~~~~~~~~~~~~~~~~~~~~~~~~~~-\frac{1}{2}\sum_{i=1}^n\sum_{a=1}^d\int_0^tD_{{h_{a,i}}}u(x_l^\PP)\beta_\PP(h_{a,i})(x_l^\PP) \dd l+M_t^u,\end{split}
\end{equation}
where $\beta_\PP(h_{a,j})$ is given in Lemma 4.5 and $M_t^u$ is a martingale with the quadratic variation process given by
\begin{equation}\label{eq4.15}\aligned\langle M_t^u\rangle&=\sum_{i=1}^n\sum_{a=1}^d\int_0^t\left\langle  X^{h_{a,i}}u(x^\PP_{r}),X^{h_{a,i}}u(x^\PP_{r})\right\rangle \dd r
\\&=\frac{1}{\varepsilon}\sum_{i=1}^n\int_0^t\left\langle \paral_{s_i}^{-1}\nabla_i f(x^\PP_{r}),\paral_{s_i}^{-1}\nabla_i f(x^\PP_{r})\right\rangle \dd r
\\&=\frac{1}{\varepsilon}\sum_{i=1}^n\int_0^t\left\langle \nabla_i f(x^\PP_{r}),\nabla_i f(x^\PP_{r})\right\rangle \dd r,\endaligned \end{equation}
with $\nabla_if(\gamma)=\nabla_i f(\gamma_{s_1},...,\gamma_{s_n})$.
\end{thm}

\beg{proof}
By \eqref{eq4.9} and applying the integration by parts formula \eqref{eq4.10} we have for $v\in \F C_0^\PP$,
\begin{equation*}
\aligned
&\E^\PP(u,v)\\&=\frac{1}{2}\sum_{i=1}^n\sum_{a=1}^d\int_{H_{\mathscr{P}}^\delta(M)} \left\langle{Du}, \hat{h}_{a,i}\right\rangle_{l^2(\mathscr{P};\mathbb{R}^d)} \left\langle{Dv}, \hat{h}_{a,i}\right\rangle_{l^2(\mathscr{P};\mathbb{R}^d)}   \dd\nu^0_\PP\\
&=\frac{1}{2}\sum_{i=1}^n\sum_{a=1}^d\int_{H_{\mathscr{P}}^\delta(M)} X^{h_{a,i}}u  X^{h_{a,i}}v \dd\nu^0_\PP\\
&=\frac{1}{2}\sum_{i=1}^n\sum_{a=1}^d\int_{H_{\mathscr{P}}^\delta(M)}v\left[ (X^{h_{a,i}}u)\beta_\PP(h_{a,i})- X^{h_{a,i}}(X^{h_{a,i}}u)\right]\dd\nu^0_\PP
\\
&=-\int_{H_{\mathscr{P}}^\delta(M)}v\left\{\frac{1}{2}\sum_{i=1}^n\sum_{a=1}^d X^{h_{a,i}}(X^{h_{a,i}}u)-
\frac{1}{2}\sum_{i=1}^n\sum_{a=1}^d\beta_\PP(h_{a,i})D_{{h_{a,i}}}u
\right\}\dd\nu^0_\PP
\\
&=-\int_{H_{\mathscr{P}}^\delta(M)}vL^\PP u\dd\nu^0_\PP,
\endaligned
\end{equation*}
where $L^\PP$ is the generator of $\E^\PP$ (see \cite[Chap. 1]{MR92}) and in the third equality we apply  \eqref{eq4.10} to $vX^{h_{a,j}}u$, which is also a smooth function on $H_\PP^\delta(M)$.
Then by the Fukushima decomposition we have under $P^z$
$$u(x_{t}^\PP)-u(x_{0}^\PP)=\int_0^t L^\PP u(x_{r}^\PP)\dd r+M_t^u,$$
where $M_t^u$ is a martingale with the quadratic variation process given by \eqref{eq4.15} (cf. \cite[Thm. 5.2.3]{FOT94}. Thus the result follows.
\end{proof}

\beg{Remark}\label{r4.9} We can also write $\E^{\mathscr{P}}$ in the following way. In fact, by \cite[Sec. 5]{AD99} we know that
for any $u=f(\gamma_{s_1},...,\gamma_{s_n})\in\F C^{\mathscr{P}}_0$,
\begin{equation}\label{eq4.16}\aligned\E^{\mathscr{P}}(u,u)=&\frac{1}{2}\int_{H_{\mathscr{P}}^\delta(M)} \langle Du,Du\rangle_{l^2(\PP,\mathbb{R}^d)} d\nu^0_\PP\\=&\sum_{i=1}^n\sum_{a=1}^d \frac{1}{2}\int_{H_{\mathscr{P}}^\delta(M)} D_{h_{a,i}}u(\gamma)D_{h_{a,i}}u(\gamma)d\nu^0_\PP\\=&\frac{1}{\varepsilon}\sum_{i=1}^n \frac{1}{2}\int_{H_{\mathscr{P}}^\delta(M)} \langle \nabla_if, \nabla_if\rangle_{T_{\gamma_{s_i}}M} d\nu^0_\PP\\=&\frac{1}{\varepsilon}\sum_{i=1}^n \frac{1}{2}\int_{M^\PP_\delta} \langle \nabla_if, \nabla_if\rangle_{T_{\gamma_{s_i}}M} d\mu_\PP,\endaligned\end{equation}
with  $$M^\PP_\delta:=\left\{\mathbf{x}=(x_{s_1},...,x_{s_n})\in M^n:\rho(x_{s_{i-1}},x_{s_i})<\delta \textrm{ for } i=1,2,...,n\right\},$$
 and $\mu_\PP:=\frac{1}{Z_\PP}\exp (-\frac{1}{2} E_\PP )\Vol_{g_\PP}$, where the energy form $E_\PP(\mathbf{x})$ is defined by
 $$E_\PP(\mathbf{x}):=\sum_{i=1}^n\frac{\rho^2(x_{s_{i-1}},x_{s_i})}{\varepsilon},$$ and
 $\Vol_{g_\PP}$ denotes the volume measure on $M^n$ with respect to the metric $g_\varepsilon^\PP:=\varepsilon g\times \varepsilon g\times ...\times \varepsilon g$.
As a result, we can also view $\E^\PP$ as a quasi-regular Dirichlet form in $L^2(M^\PP_\delta,g_\varepsilon^\PP)$.
\end{Remark}

\subsection{Derivation of the limiting process}

In order to present a better understand of the stochastic heat equation in Section 4.2. We have two ways to write the limiting equation.  The first one is invoking the stochastic parallel translation $U$:


\subsubsection{Limiting equation invoking the stochastic parallel translation $U$:}

In the following we choose $\delta>0$ satisfying the conditions in Lemma A.1 below. Then the associated finite dimension geodesic space $(H_\PP^\delta(M),G^0_\PP)$ is a smooth manifold with $nd$ dimensions.
We know that $X^{h_{a,j}}$ is a standard orthonormal frame fields in $(H_\PP^\delta(M),G^0_\PP)$ and the associated Laplace operator $\Delta_\PP$ is defined by
\begin{equation}\label{eq4.23}\Delta_\PP u=\sum_{i=1}^n\sum_{a=1}^d(X^{h_{a,i}})^2u,\quad u\in C^\infty(H_\PP^\delta(M)).\end{equation}
 Set
\begin{equation}\label{eq4.24} \beta_\PP=\frac{1}{2}\sum_{i=1}^n\sum_{a=1}^d\beta_\PP(h_{a,i})h_{a,i}.\end{equation}
Then, the generator associated to the Dirichlet form $\E^\PP$ can be written as
\begin{equation}\label{eq4.25}L_\PP=\frac{1}{2}\Delta_\PP-X^{\beta_\PP}\cdot \nabla_\PP\end{equation}
where $\nabla_\PP$ is the unique gradient associated to the metric $G^0_\PP$.
Thus,  the associated diffusion process satisfies the following equation under $P^z,$ for q.-e. $z\in H_{\PP}^\delta(M)$: for $i=1,...,n$,
\begin{equation}\label{eq4.26}\dd x_{t}^\PP(s_i)=\sum_{a=1}^dX^{h_{a,i}}(x_{t}^\PP)(s_i)\circ \dd W_t^{a,i}-X^{\beta_\PP}(x_{t}^\PP)(s_i)\dd t,\end{equation}
where $\{W^{a,i}\}$ is a sequence of independent  Brownian motions, $\circ$ means the Stratonovich integral and
\begin{equation}\label{eq4.27}\aligned X^{\beta_\PP}(\gamma)(s_i)&=\paral_{s_i}(\gamma)\beta_\PP(s_i)=\frac{1}{2}\sum_{a=1}^d\beta_\PP(h_{a,i})\paral_{s_i}(\gamma)h_{a,i}(s_i)\\
&
=\frac{1}{2\sqrt{\varepsilon}}\sum_{a=1}^d\beta_\PP(h_{a,i})\paral_{s_i}(\gamma)e_a.\endaligned\end{equation}

 As mentioned in Remark 2.5 (iii) we can also construct the $L^2$-Dirichlet form $(\E^0,\D(\E^0))$ with respect to the reference measure $\mu_0:=e^{-\frac{1}{6}\int_0^1\textrm{\Scal}(\gamma(s))\dd s}\dd \mu(\gamma).$
 We conjecture for which heuristic proofs are included in the appendix.
\vskip.10in

\noindent\textbf{Conjecture I (with limiting equation in terms of stochastic parallel translation $U$): } $(\E^\PP, \D(\E^\PP))$ Mosco converges to $(\E^0,\D(\E^0))$, as $\varepsilon\rightarrow0$. The Markov process $\Phi$ given by $(\E^0,\D(\E^0))$ satisfies the following heuristic equation
\begin{equation}\label{eq4.28}\aligned \dd  \Phi_{t,s}=
\frac{1}{2}\frac{\nabla}{ds}\partial_s \Phi_{t,s}\dd t -\frac{1}{4}\Ric(\partial_s \Phi_{t,s})\dd t+\frac{1}{12}\nabla \textrm{\Scal}(\Phi_{t,s})\dd t+ U_{t,s}(\Phi)\circ \dd W_t,\endaligned\end{equation}
where $U_{t,\cdot}$ is the stochastic parallel translation for $\Phi_t$ introduced in Section \ref{eq2.1}, $W$ is an $L^2([0,1];\mathbb{R}^d)$-cylindrical Wiener process and $\circ$ means  renormalization (see Remark \ref{r4.9}).

\subsubsection{Limiting equation invoking vector fields $\sigma$:}

Above we used the Laplace operator $\Delta_\PP$. Now we derive the diffusion equation associated to the Dirichelt form $\E^\PP$  by using the Laplacian on finite dimensional manifolds and the vector fields $\sigma$  in \eqref{eq1.1}. In this case,  $\sum_{a=1}^m\sigma_a^2$ is equal to the Laplace-Beltrami operator.  Thus by using \eqref{eq4.25} 
and \cite[Lemma 5.23]{Li}, it is easy to prove that for $u=f(\gamma_{s_1},...,\gamma_{s_n})\in\F C^{\mathscr{P}}_0$,
\begin{equation}\label{eq1}L_\PP u=\sum_{i=1}^n[\frac{1}{2\varepsilon}\Delta^{(i)}f-X^{\beta_\PP^0}(\gamma)(s_i)\cdot \nabla_if],\end{equation}
where $\Delta^{(i)}, \nabla_i$ mean the Laplace-Beltrami operator and the gradient with respect to the $i$-th variable.
\begin{equation}\label{eq 4.27}\aligned X^{\beta_\PP^0}(\gamma)(s_i)&=\frac{1}{2}\paral_{s_i}(\gamma)\frac{1}{\varepsilon} (\Delta_jb-\Delta_{j+1}b)=\frac{1}{2\varepsilon}(\gamma'(s_i-)-\gamma'(s_i+)).\endaligned\end{equation}
Therefore,   the associated diffusion process satisfies the following equation under $P^z$  q.-e. $z\in H_\PP^\delta(M)$: for $i=1,...,n$,
\begin{equation}\label{eq 4.26}\dd x_{t}^\PP(s_i)=\frac{1}{\sqrt{\varepsilon}}\sum_{a=1}^m\sigma_a(x_{t}^\PP)(s_i)\circ \dd W_t^{a,i}-X^{\beta_\PP^0}(x_{t}^\PP)(s_i)\dd t,\end{equation}
where $\{W^{a,i}\}$ is a sequence of independent  Brownian motions, $\circ$ means the Stratonovich integral.
Now we conjecture:
\vskip.10in

\noindent\textbf{Conjecture II in terms of vector fields $\sigma$: } $x^\PP$ converge to $\Phi$, as $\varepsilon\rightarrow0$, with $\Phi$ satisfying the following heuristic equation
\begin{equation}\label{eq}\aligned \dd  \Phi=
\frac{1}{2}\frac{\nabla}{ds}\partial_s \Phi\dd t + \sigma_{a}(\Phi)\circ \dd W_t^a,\endaligned\end{equation}
where  $W$ is an $L^2([0,1];\mathbb{R}^d)$-cylindrical Wiener process and $\circ$ means  renormalization (see Remark \ref{r4.9}).

\beg{Remark}\label{r4.9} $(i)$ We  only present  (so far) heuristic proofs  in the appendix. For the flat case, the convergence  can be made rigorous by classical argument (see e.g. \cite{ZZ15}). For Conjecture II we believe that the convergence above can  be made rigorous by using the theory of regularity structures introduced by Hairer in \cite{Hai14} or by using the paracontrolled distribution method proposed in \cite{GIP15}. In fact,
$$\frac{1}{2}\frac{\nabla}{\dd s}\partial_s \gamma +\frac{1}{4}\Ric(\partial_s \gamma)~~\text{and}~~U\circ \dd W, \sigma\circ \dd W$$
 are not well-defined in the classical sense and we need to multiply  two distributions. To make the proof rigorous, renormalization techniques should be involved. As there are more than 40 terms required for the renormalization for the equation \eqref{eq1.1}, the BPHZ theorem in the regularity structure theory developed  in \cite{CH16} has been used in \cite{Hai16}.  To prove the convergence rigorously in Conjecture II, the discrete version of the BPHZ theorem is required. However, there is no useful version of the discrete BPHZ theorem until now. This is one reason we do not prove  Conjecture II in this paper. We hope to be able to prove the convergence rigorously  in our future work.

 $(ii)$ We have two ways to write the limiting equations, which give us two different equations with different diffusion coefficients. Since the different approximated processes have the same law, the solutions to two different equations \eqref{eq4.28} and \eqref{eq} should have the same law.  For \eqref{eq4.28} this is more related to the integration by parts formula (see (iii) below). For \eqref{eq} this requires by regularity structure theory.

 $(iii)$
By Conjecture I,  we expect that the process given by the Dirichlet form $(\E,\D(\E))$ in Section 2.1 can be interpreted as a solution to  the following heuristic stochastic heat equation
\begin{equation}\label{eq4.37}\aligned \dd  X_{t,s}=
\frac{1}{2}\frac{\nabla}{ds}\partial_s X_{t,s}\dd t -\frac{1}{4}\Ric(\partial_s X_{t,s})\dd t+ U_{t,s}\circ \dd W_t.\endaligned\end{equation}
By the integration by parts formula by Driver in \cite{Dri92}, we can also derive (\ref{eq4.37}) heuristically: We have the following  relations: $$\frac{\nabla}{\dd s}\partial_sX_{t,s}\longleftrightarrow \int_0^1 \left\langle h'(s),\dd B_s\right\rangle,\quad\Ric(\partial_s X_{t,s})\longleftrightarrow\left\langle\int_0^1 \Ric_{U_s}h_s,\dd B_s \right\rangle,$$
for $h\in \mathbb{H}$.  Remark 4.4 gives part of the proof of the  Mosco convergence of Dirichlet forms $\E^{\mathscr{P}}$, which is equivalent to the convergence of the associated semigroups. However, the Mosco convergence in Conjecture I still requires Markov uniqueness of the limiting Dirichlet form $(\E^0, \D(\E^0))$, which is a very difficult problem in this case.

$(iii)$ If we write \eqref{eq} in local coordinates it is the same as   equation (\ref{eq1.1}) considered in \cite{Hai16}.
To use the theory of regularity structures in \cite{Hai14} or the paracontrolled distribution method in \cite{GIP15}  for equation \eqref{eq}, we may embed the manifold $M$ into a  high
dimensional Euclidean space $\mathbb{R}^N$.  In this case,  equation (\ref{eq}) can be written as
\begin{equation}\label{eq4.38}\aligned \dd  X^i=
\frac{1}{2}[\partial_{ss}^2 X^i-S^i_{jl}(X)\partial_sX^j\partial_sX^l]\dd t + \sigma^i\circ dW,\endaligned\end{equation}
where $X^i=\langle X, e_i\rangle $ with $\{e_i\}$ a basis in $\mathbb{R}^N$, $S$ is the second fundamental form and $\pi_{p}$ is the projection map from $\mathbb{R}^N$ to $T_pM$ for $p\in M$ (see also \cite{BGO}).
Here we used that
$$\frac{\nabla}{ds}\partial_s X=\partial_{ss}^2 X-S_X(\partial_sX,\partial_sX),$$
for the second fundamental form $S$ (see \cite{O83}).
By using the recent results for the theory of regularity structures in \cite{BHZ16} and \cite{CH16}, the local well-posedness of the equation (\ref{eq4.38}) follows. Moreover, by  the results in \cite{BGO} for the smooth noise case, the solution should stay in the Riemannian manifold $M$.

\end{Remark}

\beg{Remark} In \cite{AD99} another Riemannian metric $G_\PP^1$ has also been introduced and the corresponding measures  $\nu_\PP^1$ converge to the Wiener measure $\mu$. By \cite[Corollary 7.7]{AD99} we can also consider the Dirichlet form associated with $\nu_\PP^1$ and obtain that it is a quasi-regular Dirichlet form. However, it seems not easy to derive the equation for the approximation processes as in \eqref{eq4.26}.

\end{Remark}

 \appendix
  \renewcommand{\appendixname}{Appendix~\Alph{section}}
\section{Appendix}
In the appendix we give some heuristic calculations leading to proofs of  Conjectures I and II.  Before this, we prove the following results for the basis $h_{a,j}, a=1,2,...,d, j=1,...,n$.

\beg{lem}\label{l4.8} Fix  $\kappa_0, \delta$ as in Lemma 4.3 satisfying $\kappa_0\delta^2<\frac{1}{3}$ and for each $\gamma\in H^\delta_\PP(M)$, let  $b:=\phi^{-1}(\gamma)$ be the associated anti-development map and let $u:=\paral(\gamma)$ be the parallel translate of $\gamma$. Let $h:[0,\varepsilon]\rightarrow \R^d$ be the solution of the equation
\begin{equation}\label{eq4.17}h''(s)=\Omega_{u(s)}(b'(s), h(s))b'(s),\quad s\in (0,\varepsilon)\end{equation}
with boundary conditions $h(0)=0, h(\varepsilon)=\frac{1}{\sqrt{\varepsilon}}e_a,$ where $b'(s)=\frac{\Delta_ib}{\varepsilon}$ for  $\Delta_i b=b(s_i)-b(s_{i-1})$.  Then for $r\in [0,\varepsilon]$
$$h(r)=\varepsilon^{-3/2}\left(rI+\frac{[\Delta_ib]^2r^3}{6\varepsilon^2}+O(|\Delta_ib|^4)\right)\left(I-\frac{[\Delta_ib]^2}{6}+O(|\Delta_ib|^4)\right)e_a+g(r),$$
where $[\Delta_ib]^2:=\Omega_{u(0)}(\Delta_ib,\cdot)\Delta_ib$, $O(|\Delta_ib|^4)$ denotes a matrix (term) with its norm bounded by $C|\Delta_ib|^4$ and  for $r\in [0,\varepsilon]$
$$|g(r)|\leq C|\Delta_i b|^3\varepsilon^{-1/2}$$
for some constant $C$, which is independent of $\gamma$ and $r$.
\end{lem}

\beg{proof}  For convenience, let $A(s):=\Omega_{u(s)}(b'(s),\cdot)b'(s)$. By the   definition of the derivative, we know that $b'(s)=\Delta_i b/\varepsilon$.
It is easy to see that
\begin{equation}\label{eq4.18}A(u)=A(0)+\int_0^uA'(r)dr.\end{equation}
Let $\hat{h}$ be the solution to the equation \eqref{eq4.17} with $A(s)$ replaced by $A(0)$. Then it is not difficult to obtain that $\hat{h}$ satisfies the following (see \cite[Page72]{C06}):
\begin{equation}\label{eq4.19}\hat{h}(r)=\left(\sum_{n=0}^\infty\frac{A(0)^nr^{2n+1}}{(2n+1)!}\right)
\left(\sum_{n=0}^\infty\frac{A(0)^n\varepsilon^{2n+1}}{(2n+1)!}\right)^{-1}\frac{1}{\sqrt{\varepsilon}}e_a:=BD_0^{-1}\frac{1}{\sqrt{\varepsilon}}e_a.\end{equation}
Here $D_0$ is invertible, since $D_0=\varepsilon\left(I+\sum_{n=1}^\infty\frac{A(0)^n\varepsilon^{2n}}{(2n+1)!}\right):=\varepsilon(I+D)$ and
$$\|D\|\leq \sum_{n=1}^\infty \frac{\kappa_0^n|\Delta_ib|^{2n}}{(2n+1)!}\leq \frac{\kappa_0|\Delta_ib|^2}{1-\kappa_0|\Delta_ib|^2}<\frac{1}{2}\quad (\textrm{ since } \kappa_0|\Delta_ib|^2\leq \kappa_0\delta^2<1/3),$$
where we used $\|A(0)\|\leq \kappa_0|b'|^2$ in the first inequality.
Moreover, we  have
$$D_0^{-1}=\varepsilon^{-1}\sum_{n=0}^\infty (-1)^nD^n=\varepsilon^{-1}\left(I-\frac{[\Delta_ib]^2}{6}+O(|\Delta_ib|^4)\right),$$
where we used $\|\sum_{n=2}^\infty(-1)^nD^n\|\leq \sum_{n=2}^\infty\|D\|^n=\frac{\|D\|^2}{1-\|D\|}\lesssim |\Delta_i b|^4$.
In addition,
$$B=\sum_{n=0}^\infty\frac{A(0)^nr^{2n+1}}{(2n+1)!}=rI+\frac{[\Delta_ib]^2r^3}{6\varepsilon^2}+O(|\Delta_ib|^4).$$
Combining this with \eqref{eq4.19} we obtain that
$$\hat{h}(r)=\varepsilon^{-3/2}\left(rI+\frac{[\Delta_ib]^2r^3}{6\varepsilon^2}+O(|\Delta_ib|^4)\right)\left(I-\frac{[\Delta_ib]^2}{6}
+O(|\Delta_ib|^4)\right)e_a.$$
Now we give an estimate of $|h-\hat{h}|$: By Taylor's theorem we have
\begin{equation}\label{eq4.20}h(r)-\hat{h}(r)=(h'(0)-\hat{h}'(0))r+\int_0^r(r-u)[A(u)h(u)-A(0)\hat{h}(u)]\dd u,\end{equation}
which implies that
\begin{equation}\label{eq4.21}\aligned |h(r)-\hat{h}(r)|\leq& |h'(0)-\hat{h}'(0)|r+K\int_0^r(r-u)|h(u)-\hat{h}(u)|\dd u\\
&~+\int_0^r(r-u)|A(u)-A(0)||\hat{h}(u)|\dd u
\\\leq &|h'(0)-\hat{h}'(0)|r+K\int_0^r(r-u)|h(u)-\hat{h}(u)|\dd u+Cr\varepsilon^{-3/2}|\Delta_ib|^3\\:=&f(r),\endaligned\end{equation}
where in the last inequality we used \eqref{eq4.18}, that $|\hat{h}(r)|\leq C\varepsilon^{-1/2}$, and that $||A'(r)||\leq C|\Delta_ib|^3/\varepsilon^3$ from the proof of \cite[Prop. 6.2]{AD99}. Here $K:=\max ||A(r)||$. A similar argument as in the proof of Lemma \ref{L4.3} implies that
\begin{equation}\label{eq4.22}\aligned|h(r)-\hat{h}(r)|\leq &f(r)\leq \left(|h'(0)-\hat{h}'(0)|+C|\Delta_ib|^3\varepsilon^{-3/2}\right)\frac{\sinh\sqrt{K}r}{\sqrt{K}}
\\\leq&\cosh(\sqrt{\kappa_0}\delta)\varepsilon\left(|h'(0)-\hat{h}'(0)|+C|\Delta_ib|^3\varepsilon^{-3/2}\right),\endaligned\end{equation}
where we used the elementary inequality
$\frac{\sinh a}{a}\leq \cosh(a)$ for $a\in \mathbb{R}$ and
 $$ \sqrt{K}r\leq  \sqrt{\kappa_0}|\Delta_ib|r/\varepsilon\leq \sqrt{\kappa_0}|\Delta_ib|\leq\sqrt{\kappa_0}\delta .$$
Also by \eqref{eq4.20} and a similar argument as in \eqref{eq4.21} we have
$$\aligned|h'(0)-\hat{h}'(0)|\varepsilon\leq&\int_0^\varepsilon(\varepsilon-u)|A(u)h(u)-A(0)\hat{h}(u)|\dd u
\\\leq&K\int_0^\varepsilon(\varepsilon-u)|h(u)-\hat{h}(u)|\dd u+C\varepsilon^{-7/2}\int_0^\varepsilon(\varepsilon-u)u|\Delta_ib|^3\dd u
\\\leq&\cosh(\sqrt{\kappa_0}\delta)K\varepsilon\int_0^\varepsilon(\varepsilon-u)\dd u|h'(0)-\hat{h}'(0)|+CK\varepsilon^{3/2}|\Delta_ib|^3+C\varepsilon^{-1/2}|\Delta_ib|^3
\\\leq&\frac{\cosh(\sqrt{\kappa_0}\delta)\kappa_0\varepsilon|\Delta_ib|^2}{2}|h'(0)-\hat{h}'(0)|+C\varepsilon^{-1/2}|\Delta_ib|^3,\endaligned$$
where in the second inequality we used that $|\hat{h}(r)|\leq C\varepsilon^{-1/2}$, and that $||A'(r)||\leq C|\Delta_ib|^3/\varepsilon^3$. In the third inequality we used \eqref{eq4.22} and in the last inequality we used that $K\leq \kappa_0|\Delta_ib|^2/\varepsilon^2$. Since $\delta$ satisfies
$\kappa_0\cosh(\sqrt{\kappa_0}\delta)\delta^2<1$, we obtain
$$|h'(0)-\hat{h}'(0)|\leq C\varepsilon^{-3/2}|\Delta_ib|^3.$$
Therefore, combining this with \eqref{eq4.22} we have
$$|h(r)-\hat{h}(r)|\leq C\varepsilon^{-1/2}|\Delta_ib|^3,$$
which implies the result.
 \end{proof}

\noindent\textbf{Heuristic proof  of Conjectures I and II}:  We derive the convergence by heuristically analyzing the limit of the corresponding diffusion  and drift parts. The convergence of the diffusion part follows from the definition of $X^{h_{a,i}}$. In fact, for Conjecture I the diffusion part is $\paral_{s_i}(x_t^\PP)\circ \dd W^\PP$ with $W^\PP$ being an $l^2(\PP,\mathbb{R}^d)$-cylindrical Wiener process, which converges to the diffusion part heuristically. We emphasize that the diffusion part is not well-defined in the classical sense and it requires renormalization (see Remark \ref{r4.9}). For Conjecture II it is easy to see the convergence heuristically.

In the following we consider the drift part. For Conjecture I, we analyze $\frac{1}{\sqrt{\varepsilon}}\beta_\PP(h_{a,j})$:
By \eqref{eq4.10}, \eqref{eq4.6} we know that
\begin{equation}\label{eq4.29}\aligned
&\frac{1}{\sqrt{\varepsilon}}\beta_\PP(h_{a,j})(\gamma)
=\frac{1}{\varepsilon\sqrt{\varepsilon}}\left\langle \Delta_jb-\Delta_{j+1}b, h_{a,j}(s_{j}) \right\rangle +\sqrt{\varepsilon}\sum_{a_1=1}^d\left\langle
q(X^{h_{a_1,j}})h_{a_1,j},h_{a,j}\right\rangle(s_{j})
\\&=\frac{1}{\varepsilon^2}\langle \Delta_j b-\Delta_{j+1} b,e_a\rangle+\sqrt{\varepsilon}\sum_{a_1=1}^d \left\langle  h_{a,j}(s_j),  \int_{s_{j-1}}^{s_j}R_{u(r)}\big(\gamma'(r),\paral_r{h_{a_1,j}}(r)\big)\dd rh_{a_1,j}(s_j)\right\rangle\\
&:=I_1(\varepsilon)+I_2(\varepsilon),
\endaligned\end{equation}
where we used that $h_{a,j}(r)=0$ for $r\in [0,s_{j-1}]$.
Here and in the following we use $\paral_r$ to denote $\paral_r(\gamma)$ for simplicity. Now we consider $I_2(\varepsilon)$. Since
\begin{equation}\label{eq4.30}\aligned h_{a_1,j}(r)=&\frac{1}{2}h_{a_1,j}(s_{j})+\frac{1}{2}[2h_{a_1,j}(r)-h_{a_1,j}(s_{j})]:=\psi_1(r)+\psi_2(r),
\endaligned\end{equation}
 we get
\begin{equation}\label{eq4.31}\aligned I_2(\varepsilon)&=\frac{\sqrt{\varepsilon}}{2}\sum_{a_1=1}^d \left\langle  h_{a,j}(s_j),  \int_{s_{j-1}}^{s_j}R_{u(r)}\big(\gamma'(r),\paral_rh_{a_1,j}(s_j)\big)\dd rh_{a_1,j}(s_j)\right\rangle\\&~~+\sqrt{\varepsilon}\sum_{a_1=1}^d \left\langle  h_{a,j}(s_j),  \int_{s_{j-1}}^{s_j}R_{u(r)}\big(\gamma'(r),X^{\psi_2}_r\big)\dd rh_{a_1,j}(s_j)\right\rangle\\
&=:I_{21}(\varepsilon)+I_{22}(\varepsilon),\endaligned\end{equation}
where $X^{\psi_2}_r=\paral_r\psi_2(r)$.
For $I_{21}(\varepsilon)$, we  deduce that
\begin{equation}\label{eq4.32}\aligned I_{21}(\varepsilon)=&\frac{1}{2\varepsilon}\sum_{a_1=1}^d \left\langle e_a, \int_{s_{j-1}}^{s_j}\paral_r^{-1}R(\gamma'(r),\paral_re_{a_1}) \paral_re_{a_1}\dd r\right\rangle
 \\=&\frac{1}{2{\varepsilon}}\left\langle e_a, \int_{s_{j-1}}^{s_j}\paral_r^{-1}\Ric(\gamma'(r))\dd r\right\rangle
  =\frac{1}{2{\varepsilon}}\left\langle e_a, \int_{s_{j-1}}^{s_j}\Ric_{u(r)}(b'(r))\dd r\right\rangle
\\=&\frac{1}{2}\left\langle e_a, \paral_{s_j}^{-1}\Ric(\gamma'(s_j-))\right\rangle
 +\frac{1}{2{\varepsilon}}\left\langle e_a,\int_{s_{j-1}}^{s_j} \Ric_{u(r)}(b'(r))-\Ric_{u(s_j)}(b'(r))\dd r\right\rangle\\=&\frac{1}{2}\left\langle \paral_{s_j}e_a, \Ric(\gamma'(s_j-))\right\rangle
 -\frac{1}{2{\varepsilon}}\left\langle e_a,\int_{s_{j-1}}^{s_j}\int_{r}^{s_j} (D\Ric)_{u_{s}}(b'(r),b'(r))\dd s\dd r\right\rangle\\:=&I_{211}+I_{212},\endaligned\end{equation}
 where $\gamma'(s_j-)=\paral_{s_j}b'(r)$, $\Ric_{u(r)}=\paral_r^{-1}\Ric \paral_r$, $(D\Ric)_{u_{s}}(b'(s),\cdot):=(\dd /\dd s)\Ric_{u(s)}$ and  we used that $\{\paral_r e_{a_1}\}$ is an orthonormal frame in $T_{\gamma_r}M$ in the second equality. For $I_{212}$ we have
$$I_{212}=-\frac{1}{4{\varepsilon}}\left\langle e_a,(D\Ric)_{u_{s_j}}(\Delta_jb,\Delta_jb)\right\rangle+O\left(\frac{|\Delta_j b|^3}{\varepsilon}\right).$$

For $I_{22}$ we have
\begin{equation}\label{eq4.33}\aligned I_{22}(\varepsilon)=&\frac{1}{\sqrt{\varepsilon}}\sum_{a_1=1}^d \left\langle e_a, \int_{s_{j-1}}^{s_j}\paral_r^{-1}R(\gamma'(r),\paral_r\psi_2(r)) \paral_re_{a_1}\dd r\right\rangle\\=&\frac{1}{\sqrt{\varepsilon}}\sum_{a_1=1}^d \left\langle e_a, \int_{s_{j-1}}^{s_j}\Omega_{u(r)}(b'(r),\psi_2(r))e_{a_1}\dd r\right\rangle\\
 =&\frac{1}{\sqrt{\varepsilon}}\sum_{a_1=1}^d \left\langle e_a, \int_{s_{j-1}}^{s_j}\Omega_{u(s_{j-1})}(b'(r),\psi_2(r))e_{a_1}\dd r\right\rangle\\
 &+\frac{1}{\sqrt{\varepsilon}}\sum_{a_1=1}^d \left\langle e_a, \int_{s_{j-1}}^{s_j}\left(\Omega_{u(r)}(b'(r),\psi_2(r))e_{a_1}-\Omega_{u(s_{j-1})}(b'(r),\psi_2(r))e_{a_1}\right)\dd r\right\rangle\\
 :=&I_{221}(\varepsilon)+I_{222}(\varepsilon).
\endaligned\end{equation}
According to Lemma \ref{l4.8}, we know that
\begin{equation}\label{eq4.34}\aligned 2\psi_2(r)=&2h_{a_1,j}(r)-h_{a_1,j}(s_{j})=\bigg[2\varepsilon^{-\frac{3}{2}}\bigg((r-s_{j-1})I+\frac{[\Delta_jb]^2(r-s_{j-1})^3}{6\varepsilon^2}\\
 &+O(|\Delta_jb|^4)\bigg)\bigg(I-\frac{[\Delta_jb]^2}{6}+O(|\Delta_jb|^4)\bigg)
 -\frac{1}{\sqrt{\varepsilon}}\bigg]e_{a_1}+g(r-s_{j-1})\\=&:\theta(r,\varepsilon)e_{a_1}
 +O(\frac{|\Delta_jb|^2}{\varepsilon^{\frac{1}{2}}})+O(\frac{|\Delta_jb|^4}{\varepsilon^{\frac{3}{2}}}),\endaligned\end{equation}
where $\theta(r,\varepsilon)=2\varepsilon^{-\frac{3}{2}}(r-s_{j-1})-\frac{1}{\sqrt{\varepsilon}}$.
This implies that $$\aligned I_{221}(\varepsilon)=&\frac{1}{\varepsilon^{\frac{3}{2}}}\sum_{a_1=1}^d \left\langle e_a, \int_{s_{j-1}}^{s_j}\Omega_{u(s_{j-1})}(\Delta_jb,\psi_2(r))e_{a_1}\dd r\right\rangle\\
=&\frac{1}{2\varepsilon^{\frac{3}{2}}}\sum_{a_1=1}^d \left\langle e_a, \int_{s_{j-1}}^{s_j}\Omega_{u(s_{j-1})}\left(\Delta_jb,\theta(r,\varepsilon)e_{a_1}\right)e_{a_1}\dd r\right\rangle\\&+\frac{1}{\varepsilon^{\frac{3}{2}}}\bigg(O(|\Delta_jb|^3\varepsilon^{\frac{1}{2}})+O(|\Delta_jb|^5\varepsilon^{-\frac{1}{2}})\bigg)
\\=&O(|\Delta_jb|^3\varepsilon^{-1})+O(|\Delta_jb|^5\varepsilon^{-2}),\endaligned$$
where we used that $\int_{s_{j-1}}^{s_j}\theta(r,\varepsilon)\dd r=0$ in the third equality.
 In the following, we use \eqref{eq4.34} to estimate $I_{222}(\varepsilon)$:
\begin{equation}\label{eq4.35}\aligned I_{222}(\varepsilon)=&
 \frac{1}{\sqrt{\varepsilon}}\sum_{a_1=1}^d \left\langle e_a, \int_{s_{j-1}}^{s_j}\left(\Omega_{u(r)}(b'(r),\psi_2(r))e_{a_1}-\Omega_{u(s_{j-1})}(b'(r),\psi_2(r))e_{a_1}\right)\dd r\right\rangle
 \\=&
 \frac{1}{2\sqrt{\varepsilon}}\left\langle e_a, \int_{s_{j-1}}^{s_j}\theta(r,\varepsilon)\left[\Ric_{u(r)}(b'(r))-\Ric_{u(s_{j-1})}(b'(r))\right]\dd r\right\rangle+O_{b,\varepsilon}
\\ =&
 \frac{1}{2{\sqrt{\varepsilon}}}\left\langle e_a, \int_{s_{j-1}}^{s_j}\theta(r,\varepsilon)\int_{s_{j-1}}^{r}\left(D\Ric\right)_{u(s)}(b'(s),b'(s))\dd s\dd r\right\rangle+O_{b,\varepsilon}\\
 =&\frac{1}{2{\varepsilon^{\frac{5}{2}}}}\left\langle e_a,
 \left(D\Ric\right)_{u_{s_j}}(\Delta_jb,\Delta_jb)\right\rangle \int_{s_{j-1}}^{s_j}\left[\frac{2(r-s_{j-1})}{\varepsilon^{3/2}}-\frac{1}{\sqrt{\varepsilon}}\right](r-s_{j-1})\dd r+O_{b,\varepsilon}
 \\=&\frac{1}{12\varepsilon} \left\langle e_a, \left(D\Ric\right)_{u_{s_j}}(\Delta_jb,\Delta_jb)\right\rangle+O_{b,\varepsilon},
\endaligned\end{equation}
with $O_{b,\varepsilon}:=O\left(\frac{|\Delta_jb|^3}{\varepsilon}\right)+O\left(\frac{|\Delta_j b|^5}{\varepsilon^2}\right)$.
Combining the computation for $I_{212}$ and $I_{221}, I_{222}$ we have
$$I_{212}+I_{22}=-\frac{1}{6\varepsilon} \left\langle e_a, \left(D\Ric\right)_{u_{s_j}}(\Delta_jb,\Delta_jb)\right\rangle+O_{b,\varepsilon},$$
which combined with \eqref{eq4.29}, \eqref{eq4.31}, \eqref{eq4.32} implies that
$$\aligned
\frac{1}{\sqrt{\varepsilon}}\beta_\PP(h_{a,j})(\gamma)
=&\frac{1}{\varepsilon^2}\langle \Delta_j b-\Delta_{j+1} b,e_a\rangle+\frac{1}{2}\left\langle \paral_{s_j}e_a, \Ric(\gamma'(s_j-))\right\rangle\\&-\frac{1}{6\varepsilon} \langle(D\Ric)_{u_{s_j}}(\Delta_jb,\Delta_jb),e_a\rangle+O_{b,\varepsilon}.
\endaligned$$ Then we obtain
$$\aligned X^{\beta_\PP}(\gamma)(s_i)=&\frac{1}{2\varepsilon^2}\sum_{a=1}^d\langle \Delta_j b-\Delta_{j+1} b,e_a\rangle\paral_{s_j}e_a+\frac{1}{4}\sum_{a=1}^d\left\langle \paral_{s_j}e_a, \Ric(\gamma'(s_j-))\right\rangle\paral_{s_j}e_a\\&-\frac{1}{12\varepsilon} \sum_{a=1}^d\langle(D\Ric)_{u_{s_j}}(\Delta_jb,\Delta_jb),e_a\rangle \paral_{s_j}e_a+O_{b,\varepsilon}\\=&\frac{1}{2\varepsilon^2}\paral_{s_j}( \Delta_j b-\Delta_{j+1} b)+\frac{1}{4}  \Ric(\gamma'(s_j-))-\frac{1}{12\varepsilon} \paral_{s_j}(D\Ric)_{u_{s_j}}(\Delta_jb,\Delta_jb)+O_{b,\varepsilon},\endaligned$$
where we used that $\{\paral_{s_j}e_a\}$ is an orthonormal frame  in $T_{\gamma_{s_j}}M$.
Heuristically, we have $\frac{1}{\varepsilon^2}\paral_{s_j}( \Delta_j b-\Delta_{j+1} b)\rightarrow \paral_s\partial_{ss}^2b$ and $O_{b,\varepsilon}\rightarrow0$ as $\varepsilon\rightarrow0$, since $|\Delta_j b|\simeq \varepsilon^{\frac{1}{2}-}$. We also have $\sum_{i=1}^d\paral_{s_j}(D\Ric)_{u_{s_j}}(e_i,e_i)=\nabla \textrm{Scal}$, which suggests the third term above converges to $\frac{1}{12}\nabla \textrm{\Scal}(\gamma)$ by a similar argument as in \cite[Section 6]{AD99}. Now we have
\begin{equation}\label{eq4.36}\aligned X^{\beta_\PP}(\gamma)\rightarrow&-\frac{1}{2}\paral_{s}\partial_{ss}^2b +\frac{1}{4}\Ric(\partial_s \gamma)-\frac{1}{12}\nabla \textrm{\Scal}(\gamma)\\=&-\frac{1}{2}\frac{\nabla}{ds}\partial_s \gamma +\frac{1}{4}\Ric(\partial_s \gamma)-\frac{1}{12}\nabla \textrm{\Scal}(\gamma),\endaligned\end{equation}
where we used $\paral_{s}\partial_{ss}^2b=\paral_{s}\frac{\dd }{\dd s}\paral_{s}^{-1}\partial_s\gamma=\frac{\nabla}{\dd s}\partial_s \gamma$.

For Conjecture II it is easy to see that $\frac{1}{\varepsilon^2}\paral_{s_j}( \Delta_j b-\Delta_{j+1} b)\rightarrow \paral_s\partial_{ss}^2b=\frac{\nabla}{\dd s}\partial_s \gamma$.

\section*{Acknowledgments}
We are very grateful to Professor Fengyu Wang and Dr. Xin Chen for numerous discussions which helped us understand  concepts in stochastic analysis on Riemannian manifolds. We also would like to thank  Dr. Zhehua Li for pointing out \cite[Lemma 5.23]{Li} and \eqref{eq1} to us. We would also like to thank Professor Martin Hairer and Professor Xuemei Li for inspiring discussions.

\beg{thebibliography}{99}

\bibitem{A96} S. Aida, \emph{Logarithmic Sobolev inequalities on loop spaces over compact Riemannian
manifolds}, in ``Proc. Fifth Gregynog Symp., Stoch. Anal. Appl. ''(I. M. Davies, A. Truman,and K. D. Elworthy, Eds.), pp. 1-15, World Scientific, 1996

\bibitem{AD99}  L. Andersson,  B. K. Driver, \emph{Finite-dimensional approximations to
Wiener measure and path integral formulas on manifolds},
J. Funct. Anal, 165, no. 2, (1999),430-498.

\bibitem{AE95}S. Aida and D. Elworthy, \emph{Differential calculus on Path and Loop spaces 1: Logarithmic
Sobolev inequalities on Path spaces,} C.R. Acad. Sci. Paris Ser. I Math. 321 (1995),
97-102.

\bibitem{ALR93}S. Albeverio, R. L\'{e}andre and M. R\"{o}ckner, \emph{Construction of rotational invariant diffusion
on the free loop space,} C. R. Acad. Paris Ser. 1 316 (1993), 287-292

\bibitem{AR91}S. Albeverio, M. R\"{o}ckner, \emph{Stochastic differential equations in infinite
dimensions: Solutions via Dirichlet forms}, Probab. Theory Related Field 89 (1991), 347-386

\bibitem{BE85} D. Bakry, M. \'{E}mery, Diffusions hypercontractives, in S\'{e}minaire de probabilit\'{e}s, XIX,
1983/84, vol. 1123 of Lecture Notes in Math., Springer, Berlin, (1985), 177-206

\bibitem{BGO} Z. Brze\'{z}niak, Goldys B.and M. Ondrej\'{a}t: Stochastic geometric partial differential equations. New trends
in stochastic analysis and related topics, 1-32, Interdiscip. Math. Sci., 12, World Sci. Publ., Hackensack,
NJ, 2012

\bibitem{BHZ16}Y. Bruned, M. Hairer, and L. Zambotti. \emph{Algebraic renormalisation of regularity structures}
{\it arXiv:1610.08468,}

\bibitem{CG02}A. Chojnowska-Michalik, B. Goldys, Symmetric Ornstein-Uhlenbeck semigroups and their generators, Probab. Th. Relat. Fields, 124(2002), 459-486

\bibitem{CHL} B. Capitaine,  E. P. Hsu and M. Ledoux, \emph{Martingale representation and a simple proof of logarithmic
Sobolev inequalities on path spaces,} Electron. Comm. Probab. 2(1997), 71-81.

\bibitem{CH16}A. Chandra and M. Hairer. \emph{An analytic BPHZ theorem for regularity structures}
{\it arXiv:1612.08138}

\bibitem{C06}I. Chavel, Riemannian geometry- A modern introduction, in "Cambridge Tracts in
Mathematics," Vol. 108, Cambridge Press, Cambridge, second edition, 2006.

\bibitem{CLW11}X. Chen, X.-M. Li and B. Wu, \emph{A concrete estimate for the weak Poincar´e inequality on
loop space,} Probab.Theory Relat. Fields 151 (2011), no.3-4, 559-590.

\bibitem{CLW17}X. Chen, X.-M. Li and B. Wu, \emph{Functional inequalities on general loop spaces,} Preprint.

\bibitem{CLWFL}X. Chen, X.-M. Li and B. Wu, \emph{Analysis on Free Riemannian Loop Space,} Preprint.

\bibitem{CW14} X. Chen, B. Wu,  \emph{Functional inequality on path space over a non-compact Riemannian manifold,} J. Funct. Anal. 266 (2014), 6753-6779

\bibitem{CT} L.J. Cheng, A. Thalmaier, \emph{Characterization of pinched Ricci curvature by functional inequalities,} J. Geom. Anal. 28(2018), 2312-2345

\bibitem{D04} G. Da Prato, Kolmogorov Equations for Stochastic PDEs, Birkh\"{a}user, 2004.

\bibitem{DR92}B. K. Driver and M. R\"{o}ckner,\emph{Construction of diffusions on path and loop spaces of
compact Riemannian manifolds,} C. R. Acad. Sci. Paris Ser. I 315 (1992), 603-608

\bibitem{Dri92}B. K. Driver, \emph{A Cameron-Martin quasi-invariance theorem for Brownian motion on
a compact Riemannian manifolds,} J. Funct. Anal. 110 (1992), 273-376.

\bibitem{Dri94} B. K. Driver, \emph{A Cameron-Martin type quasi-invariance theorem for
pinned Brownian motion on a compact Riemannian manifold}, Trans. Amer.
Math. Soc. 342, no. 1, (1994), 375-395.

\bibitem{EL99} K. D. Elworthy, X.- M. Li and Y. Lejan, \emph{On The geometry of diffusion
operators and stochastic flows,} Lecture Notes in Mathematics, 1720(1999), Springer-Verlag.

\bibitem{EM} K. D. Elworthy and Z.-M. Ma, \emph{Vector fields on mapping spaces and
related Dirichlet forms and diffusions,} Osaka. J. Math. 34(1997), 629-651.

\bibitem{F94} S.- Z. Fang, \emph{Un in\'equalit\'e du type Poincar\'esur un espace de chemins,}
C. R. Acad. Sci. Paris S\'erie I 318(1994), 257-260.

\bibitem{Fun88} T. Funaki, \emph{On diffusive motion of closed curves,} In Probability theory
and mathematical statistics (Kyoto, 1986), vol. 1299 of Lecture Notes
in Math., 86-94. Springer, Berlin, 1988.

\bibitem{Fun92} T. Funaki, \emph{A stochastic partial differential equation with values in a manifold},
J. Funct. Anal. 109(1992), 257-288.

\bibitem{FOT94}M. Fukushima, Y. Oshima, and M. Takeda, \emph{Dirichlet Forms and Symmetric
Markov Processes,} de Gruyter, Berlin (1994), second edition, (2011)

\bibitem{FW05}S. Fang and F.-Y. Wang, \emph{Analysis on free Riemannian path spaces,} Bull. Sci. Math.
129 (2005), 339-355.

 \bibitem{FW15}S. Z. Fang, B. Wu,  \emph{Remarks on spectral gaps on the Riemannian path space,} Elect. Commun. Probab.,
 22 (2017), 1-13

\bibitem{GIP15} M. Gubinelli, P. Imkeller, N. Perkowski, \emph{Paracontrolled distributions and singular PDEs}, Forum Math. , 6(2015), 1-75

\bibitem{GM98} F. Z. Gong, Z. M. Ma, \emph{The log-Sobolev inequality on loop space over a
compact Riemannian manifold}, J. Funct. Anal. 157 (1998), 599-623.

\bibitem{GRW01} F. Z. Gong, M. R\"{o}ckner, L. Wu \emph{Poincar\'{e} Inequality for Weighted First Order Sobolev
Spaces on Loop Spaces}, J. Funct. Anal. 185(2001), 527-563

\bibitem{GW06} M. Gourcy,  L. Wu, \emph{Logarithmic Sobolev inequalities of diffusions for the $L^2$-metric}. Potential
 Anal. 25(2006), 77-102

\bibitem{Hai14}M. Hairer. \emph{A theory of regularity structures} Invent. Math., 198(2):269-504, 2014

\bibitem{Hai16}M. Hairer. \emph{The motion of a random string,} arXiv:1605.02192

\bibitem{H95} E. P. Hsu, \emph{Quasi-invariance of the Wiener measure on the path space over a compact
Riemannian manifold,} J. Funct. Anal. 134(1995), 417-450.

\bibitem{H02} E. P. Hsu, \emph{Stochastic Analysis on Manifold,} American Mathematical Society, 2002.

\bibitem{H97}E.P. Hsu, \emph{Logarithmic Sobolev inequalities on path spaces over Riemannian manifolds,} Comm. Math.
Phys. 189 (1997) 9–16.

\bibitem{H99} E. P. Hsu, \emph{Analysis on path and loop spaces,} in "Probability Theory and Applications"
(E. P. Hsu and S. R. S. Varadhan, Eds.), LAS/PARK CITY Mathematics Series, 6(1999), 279--347,
Amer. Math. Soc. Providence.

\bibitem{IM85} A. Inoue, Y. Maeda, \emph{On integral transformations associated with a
certain Lagrangian—as a prototype of quantization,} J. Math. Soc. Japan 37,
no. 2, (1985), 219-244

\bibitem{K06}A. V. Kolesnikov. \emph{Mosco convergence of Dirichlet forms in infinite dimensions with changing
reference measures}. J. Funct. Anal., 230(2006), 382-418

\bibitem{Li}Z.H., Li. \emph{A finite dimensional approximation to pinned Wiener measure on some symmetric spaces,} {\it arXiv: 1306.6512v4}.

\bibitem{L} J.-U. L\"{o}bus, \emph{A class of processes on the path space
over a compact Riemannian manifold with unbounded diffusion,} Tran.
Ame. Math. Soc. (2004), 1-17.

\bibitem{MR92}Z. M. Ma and M. R\"{o}ckner, \emph{Introduction to the Theory of (Non-Symmetric) Dirichlet
Forms ,}(Springer-Verlag, Berlin, Heidelberg, New York, 1992.

\bibitem{MR00}Z. M. Ma and M. R\"{o}ckner, \emph{Construction of diffusions on configuration spaces,} Osaka
J. Math. 37 (2000) 273-314.

\bibitem{N} A. Naber, \emph{Characterizations of bounded Ricci curvature on smooth and nonsmooth spaces}, {\it arXiv: 1306.6512v4}.

\bibitem{Nor98} J. R. Norris,\emph{Ornstein-Uhlenbeck processes indexed by the circle}, Ann.
Probab. 26, no. 2, (1998), 465-478.

\bibitem{Nor95} J. R. Norris,\emph{Twisted sheets}, J. Funct. Ana. 132(1995), 273-334.

\bibitem{O83}  B. O'Neill, Semi-Riemannian Geometry. with application to relativity, Pure and
Applied Mathematics 103, Academic Press, Inc., New York, 1983.

\bibitem{RZZ15}M. R\"{o}ckner, R. Zhu, X. Zhu, \emph{Restricted Markov unqiueness for the
stochastic quantization of $P(\phi)_2$ and its
applications}, J. Funct. Anal.  272 (2017), 4263-4303

\bibitem{RZZ16}M. R\"{o}ckner, R. Zhu, X. Zhu, \emph{Ergodicity for the stochastic quantization problems on the 2D-torus},  Commun. Math. Phys. 352, 1061–1090 (2017)

\bibitem{TL} T. Laetsch, \emph{An approximation to Wiener measure and quantization of the Hamiltonian on manifolds with non-positive sectional curvature,} J. Funct. Anal. 265 (2013), no. 8, 1667–1727. MR 3079232 1.1, A

\bibitem{W04} F.- Y. Wang, \emph{Weak poincar\'{e} Inequalities on path
spaces,} Int. Math. Res. Not. 2004(2004), 90--108.

\bibitem{W05} F.Y. Wang, \emph{Functional Inequalities, Markov Semigroup and Spectral Theory}. Chinese Sciences
Press, Beijing (2005)

\bibitem{W14}  F.- Y. Wang, \emph{Analysis for diffusion processes on Riemannian manifolds,} World
Scientific, 2014.

\bibitem{W17}  F.- Y. Wang, \emph{Identifying constant curvature manifolds, Einstein manifolds, and Ricci parallel manifolds,} {\it arXiv:1710.00276v2}

\bibitem{WW08}F.-Y. Wang and B. Wu, \emph{Quasi-regular Dirichlet forms on path and loop spaces,} Forum
Math. 20 (2008), 1085-1096.

\bibitem{WW09} F. -Y. Wang and B. Wu, \emph{Quasi-Regular Dirichlet Forms on Free Riemannian Path and Loop Spaces,} Inf. Dimen. Anal. Quantum Probab.
and Rel. Topics 2(2009), 251-267.

\bibitem{WW16}F.- Y. Wang and B. Wu, \emph{Pointwise Characterizations of Curvature and Second
Fundamental Form on Riemannian Manifolds,} {\it arXiv:1605.02447.}

\bibitem{W1} B. Wu, \emph{Characterizations of the upper bound of Bakry-Emery curvature,} {\it arXiv:1612.03714}

\bibitem{W2} B. Wu, \emph{Finite Dimensional Approximations to Wiener Measure on General Noncompact Riemannian Manifold,} preprint.
\bibitem{ZZ15} R. Zhu, X. Zhu, \emph{Lattice approximation to the dynamical $\Phi_3^4$ model},  The Annals of Probability, 46 (2018), no. 1, 397-455.

\end{thebibliography}
\end{document}